\documentclass[11pt]{amsart}
\usepackage[utf8]{inputenc}
\usepackage{geometry}                
\usepackage{graphicx}
\usepackage{amssymb}
\usepackage{epstopdf}
\usepackage{fancyhdr}
\usepackage[colorlinks=true, pdfstartview=FitV, linkcolor=blue, 
            citecolor=blue, urlcolor=blue]{hyperref}

\newtheorem{theorem}{Theorem}[section]

\newtheorem{lemma}{Lemma}[section]

\newtheorem{prop}{Proposition}[section]

\newcommand{\pproof}[2]{{\raggedleft \bf {Proof #1}:}\\#2  \begin{flushright} $\square$ \end{flushright}}

\newcommand{\I }{\mathcal{I}}
\newcommand{\J }{\mathcal{J}}
\newcommand{\vv }{\upsilon}
\newenvironment{dema}{\ \\ {\bf Proof}}{\hfill\mbox{$\square$}\medskip}

\newcommand{\R }{\mathbb{R}}
\newcommand{\PP }{\mathbb{P}}
\newcommand{\E }{\mathbb{E}}
\newcommand{\N }{\mathbb{N}}

\newcommand{\F }{\mathcal{F}}

\newcommand{\T }{\mathcal{T}}

\newcommand{\BS }{Black $\&$ Scholes}

\newcommand{\be}{\begin{equation*}}
\newcommand{\ben}{\begin{equation}}
\newcommand{\ee}{\end{equation*}}
\newcommand{\een}{\end{equation}}

\newtheorem{remark}{Remark}

\title{The  critical price of the American put near maturity in the Jump Diffusion model}

\begin{document}
\maketitle
\vspace {-0.cm}
\begin{center}
\renewcommand{\thefootnote}{(\arabic{footnote})}
  \scshape Aych Bouselmi\footnote{{aych.bouselmi@gmail.fr}}, Damien Lamberton \footnote{damien.lamberton@univ-mlv.fr}
\renewcommand{\thefootnote}{\arabic{footnote}}\setcounter{footnote}{0}
\end{center}
\vglue0.3cm
\hglue0.02\linewidth
\begin{minipage}{0.9\linewidth}
\begin{center}
Universit\'e Paris-Est\\
Laboratoire  d'Analyse et de Mathématiques Appliquées (UMR 8050), UPEM\\
UPEC, CNRS, Projet Mathrisk INRIA,\\
F-77454, Marne-la-vallée, FRANCE\\
\end{center}
\end{minipage}
\vspace{1cm}
\begin{abstract}
We study the behavior of the critical price of an American put option near maturity in  the Jump diffusion  model when the underlying stock pays dividends at a continuous rate and  the limit of the critical price is smaller than the stock price. In particular, we prove that, unlike the case where the limit is equal to the strike price, jumps can  influence  the convergence rate.
   \end{abstract}

\section*{Introduction}
The behavior of the critical price of the American put near maturity has been deeply investigated. Its limit  was characterized  in the Black Scholes model  (see  \cite{kim,moerbeke}) by $$b(T):= \lim_{t\to T}b(t)=
  \min \left({r\over\delta}K, K\right),$$ where $r$ and $\delta$ denote the interest rate and the dividend rate
  and $b(t)$ is the critical price at time $t$.

 This result was  generalized to more general  exponential Lévy models in \cite{papiermikou}. In fact, denoting $\bar d= r-\delta -\int (e^y-1)^+\nu(dy)\;\footnote{The quantity $\bar d$ is denoted by $ d^+$ in \cite{papiermikou}}$, with $\nu$  the Lévy measure of the underlying Lévy process, we have 
 $$b(T)=K\:,\; \mbox{if} \quad\bar d\geq 0, $$
 and  $$b(T)=\xi\;,\; \mbox{if}\quad \bar d< 0, $$
  where  $\xi$ is the unique  solution, in $[0,K]$, of
\begin{equation}\label{eq.carac}
rK-\delta x-\int (xe^y-K)^+\nu(dy)=0.
\end{equation}

In the Black Scholes Model, the quantity $\bar d$  reduces to $\bar d=r-\delta$ and we  distinguish, according as $\bar d>0$, $\bar d=0$ and $\bar d<0$,  different behaviors of the critical price near maturity. 
 In fact, Barles et al in \cite{barles} (see also D. Lamberton \cite{damien}) established, in the case where $\bar d>0$ (which implies $b(T)=K$), that
\begin{equation}\label{95's}
\frac{K-b(t)}{\sigma K}\sim_{t\rightarrow T}  \sqrt{(T-t)|\ln(T-t)|},
\end{equation}
where the expression  $f\sim_{t\to a}g$ (or $f\sim_{ a}g$) is equivalent to $\lim_{t\to a}{f(t)\over g(t)}=1$.
The cases $\bar d<0$ and $\bar d=0$ were investigated by D. Lamberton and S. Villeneuve in \cite{villeneuve} and they obtained
: \\
If  $\bar d=0$ (which also implies $b(T)=K$)
$$\frac{K-b(t)}{\sigma K}\sim_{t\rightarrow T}  \sqrt{2(T-t)|\ln(T-t)|}.$$
If  $\bar d<0$ ($b(T)<K$), there exists $y_{0}\in (0,1)$, which is characterized thanks to an auxiliary optimal stopping problem, such that
$$\frac{b(T)-b(t)}{\sigma b(T)}\sim_{t\rightarrow T} y_{0} \sqrt{(T-t)}.$$
The critical price has  also been studied in the Jump diffusion model. In fact, Pham proved in \cite{pham} 
that the result \eqref{95's}, obtained in \cite{barles,damien},
remains exactly the same in the Jump diffusion model, in the case where $\bar d>0$ and $\delta=0$.
This remains true if $\delta>0$ (see \cite{mikou}). 

The purpose of this paper is to study the convergence rate of the critical  price of the American put, in the Jump diffusion model,  with $\bar d\leq0$.
Considering the results of Pham in \cite{pham}, we expect to obtain the same results as the study performed  by Lamberton and Villeneuve in  the Black-Scholes model when ($\bar d=r-\delta\leq0$ ), meaning that jumps do not have any  influence on the convergence rate.
Surprisingly, we obtain the expected result only  for the case $\bar d=0$. Indeed, we obtain for $\bar d=0$ (see Theorem \ref{thjd6}),
$$\frac{K-b(t)}{\sigma K}\sim_{t\rightarrow T}  \sqrt{2(T-t)|\ln(T-t)|},$$
and for $\bar d<0$ (see Theorem \ref{thjd2}),
$$\frac{b(T)-b(t)}{\sigma b(T)}\sim_{t\rightarrow T} y_{\lambda,\beta}\sqrt{(T-t)},$$
where $y_{\lambda,\beta}$ is a real  umber satisfying $y_{\lambda,\beta}\geq y_{0}$, 
and depending on $\nu(\{\frac{K}{b(T)}\})$ we can have $y_{\lambda,\beta} > y_{0}$. This point will be discussed in more details in section \ref{S2}.

This study is composed of four sections. 
In Section 1, we recall  some useful  results on the American put  which will be used throughout this study. 
In Section 2, we  give some results on the regularity   of the American put price and the early  exercise premium.
In Section 3, we investigate the case where the limit of the critical price is far from the singularity $K$. 
Therefore,   we  have enough regularity to give an expansion of the American put price near maturity from which
  the  critical price behavior  will be deduced.  The method is similar to the one used in \cite {villeneuve}
  and is based on an expansion of the American put price along parabolas. 
  However, the possibility that the stock price jumps into a neighborhood of the exercise price produces 
  a contribution of the local time in the expansion.
Section 4 is devoted to the study of the case $\bar d=0$. In this case $b(T)=K$, hence we have no longer enough smoothness to  obtain an expansion around the limit point $(T,b(T))$. Then we will study the behavior of the European critical price $b_e(t)$ instead of $b(t)$. Thereafter, we prove that $b(t)$ and $b_e(t)$ have the same behavior.

\section{Preliminary }

In the Jump Diffusion model, under a risk-neutral probability, the risky asset price is modelized by   
$\left (S_t\right )_{t\geq0}$ given by
$$
S_t=S_0e^{\tilde{X}_t},\qquad\mbox{ with }\tilde{X}_t=(r-\delta)t+\sigma B_t-\frac{\sigma^2}{2}t+Z_t-t\int (e^y-1)\nu(dy)
$$
where $r>0$ is the interest rate, $\delta\geq0$ the dividend rate, $\left (B_t\right )_{t\geq0}$ a Standard Brownian Motion  and 
 $\left (Z_t\right )_{t\geq0}$ a Compound Poisson Process and $\nu$ its Levy measure.
We then have
$$
dS_t   =S_{t^-}\left(\gamma_0dt+\sigma dB_t+d\bar Z_t\right),
\mbox{ with }
\bar Z_t=\Sigma_{0<s\leq t}(e^{\Delta Z_s}-1)
\mbox{ and }
\gamma_0=r-\delta-\int(e^y-1)\nu(dy).
$$
Denote by $\mathbb{F}$ the completed natural filtration  of the process $\tilde X_t$
and  suppose all over this paper  that  the following  assumptions are satisfied
 $$\sigma>0, \qquad \nu(\R)<\infty,\qquad\int e^y\nu(dy)<\infty
\qquad\mbox{and}\qquad \bar d=r-\delta-\int_{y>0}(e^y-1)\nu(dy)\leq0.$$
 The price of an American put with maturity $T>0$ and strike price  $K>0$ is given, at  $t\in[0,T]$, by $P(t,S_t)$ with $P$ defined  for all $(t,x)\in [0,T]\times\R^+$ by
$$P(t,{x})=\sup_{\tau\in \mathcal{T}_{0,T-t}}\E(e^{-r\tau}(K-xe^{\tilde X_{\tau}})_+),$$
where $\mathcal{T}_{0,T-t}$ is the set of all $\mathbb{F}$-stopping times  taking values in $[0,T-t]$.
 The value function $P$ can also be characterized (see \cite{papiermikou}) as the unique continuous and bounded solution of the following variational inequality
$$\max\{\psi-P;\frac{\partial P}{\partial t}+\mathcal{A} P -rP\}= 0,\;  \mbox{ (in the  sense of  distributions)},$$
 with the terminal condition $P(T,.)=\psi$. Here $\mathcal{A}$ is the infinitesimal generator of the process $S$.
   The free boundary of this variational inequality is called the exercise boundary,
and  at each $t\in[0,T]$, the critical price is given by $$b(t)=\inf\left \{x>0 \;|\; P(t,x)>(K-x)^+  \right \}.$$ It was proved in \cite{papiermikou} that, if $\bar d\leq 0$, then
\begin{equation}
\lim_{t\rightarrow T}b(t)=\xi:=b(T),
\end{equation}
  where $\xi$ is the unique solution, in $[0,K]$, of $rK=\delta x+\int(xe^y-K)^+\nu(dy)$. Note that, if $\bar d=0$, then $b(T)=\xi=K$.

Finally, recall that the price of a European put with   maturity $T$ and   strike price $K$  is given, at time $t$, by 
$$P_e(t,x)=\E\left(e^{-r(T-t)}(K-S_{T-t})_+\;|\; S_0=x\right).$$
The quantity $(P-P_e)$ is called the early exercise premium, we then have $P(t,x)=P_e(t,x)+e(T-t,x).$ 
 Setting $\theta=T-t$, then   the early  exercise prime, $e(\theta,x)$, is characterized for the American put in the exponential Levy model as follows (see\cite{mikou})
\begin{eqnarray*}
&&e(\theta ,x)=\\
&&\E \left\{
\int_0^\theta e^{-rs}
\left(
 rK-\delta S^x_s-\int_{y>0}\left [ P(t+s,S_s^xe^y)-\left(K-S_s^x e^y\right)\right ]\nu(dy) 
\right) 
 1_{\{S_s^x<b(t+s)\}} ds
 \right\}.
\end{eqnarray*}
We also define, for all $t\in \left (0, T\right )$, the  European critical  price,  $b_e(t)$,  as the unique solution of 
\[F(t,x)= P_e(t,x)-(K-x)=0.\]
It  easy to check that, for all $t\in \left (0, T\right )$, $b_e(t)$ is well defined,  $b_e(t)\in\left (0, K\right )$.  It is also straightforward that $P_e\leq P$,  therefore $b(t)\leq b_e(t)\leq K$ .

\section{Regularity estimate for the value function in the jump diffusion model}

In this section, we study the spatial derivatives behavior of $P$, $P_e$ and $e(\theta,x)$ near $(T, b(T))$. 
We also give a lower bound for the second spatial derivative near $(T,b(T))$. These results will be proved in Appendix 1.
\begin{lemma}\label{deriv_prime}
Under the model assumption, we have 
\begin{enumerate}
\item For all $x\in  \left(0,b_e(t)\wedge b(T)\right]$, we have,
as $\theta(=T-t)$ goes to 0,
 $$\left |\frac{\partial e}{\partial x} (\theta,x)\right |=\frac{1}{x} o(\sqrt{\theta}),$$ with $o(\sqrt\theta)$ 
 uniform with respect to $x$.
\item For all  $x\in \left(0,  b(T)\wedge b_e(t)\right]$, we have $$\frac{\partial P}{\partial x}(t,x)+1=(1+\frac{1}{x})o(\sqrt{\theta}),$$ with $o(\sqrt\theta)$ uniform with respect to $x$.
\end{enumerate}
\end {lemma}
\begin{lemma}\label{min_deriv2_P}
According to the hypothesis of the model,  we have, for all $b(t) \leq x < b(T)\wedge b_e(t)$ and for all  $\theta=T-t$ small enough, the following inequality
\begin{eqnarray*}
\inf_{b(t)< u < x} \frac{u^2\sigma^2}{2}\frac{\partial^2P}{\partial x^2}(t,u)\geq
&& \left (\bar\delta-\epsilon (\theta)\right )\left (b(T)-x \right ) -{\lambda\beta} \E\left (\sigma B_\theta-\ln(\frac{b(T)}{x})\right )^++o(\sqrt{\theta}),\\
\end{eqnarray*}
with   $\lim_{\theta \downarrow 0}\epsilon(\theta)=0$,~~ $\bar\delta=\delta +\int_{\{y>\ln(\frac{K}{b(T)}\}} e^y\nu(dy),\; {\lambda}=\nu\left\{\ln\left(\frac{K}{b(T)} \right)\right\}$  ~and~${\beta}=\frac{K}{\bar\delta b(T)}.$
\end{lemma}
\section{Regular case}
We begin this section with introducing an auxiliary optimal stopping problem which will be needful for deriving the expansion of the American put price near maturity along a parabolic branch. Once we have this expansion we will be able to derive the convergence rate of the critical price.
\subsection{{An auxiliary optimal stopping problem}}\label{etud_u_alpha}
Let  $ {\beta} $  be a non-negative number,  $(B_s)_{s\geq 0}$ be  a standard Brownian motion with  local time at $x$ 
denoted by $\tilde L^x$. We denote by $\T_{0,1}$ the set of all $\sigma \left (B_t \;;\;t\geq0\right )$-stopping times with
 values in $[0,1]$ . Consider also a Poisson process $(N_s)_{s\geq 0}$,   independent of $B$,  with intensity $\lambda$, 
 we denote by $\hat T_1$ its first jump time and by $\hat \T_{0,1}$ the set of all 
 $\sigma \left ((N_t,B_t) \;;\;t\geq0\right )$-stopping times with values in $[0,1]$.  
 We define the  functions $\vv_{\lambda, \beta}$ 
as follows
\begin{eqnarray*}
{v}_{\lambda,\beta}(y)&=&
 \sup_{\tau\in \hat\T_{0,1}}
\E\left[ e^{\lambda\tau} 1_{\{\hat N_\tau=0\}}
\int_0^{\tau} f_{\lambda\beta}(
   y+ B_{s}
    )ds
     +\frac{\beta}{2} 
     e^{\lambda\tau} 1_{\{\hat N_\tau=1\}}
      \left(L^{-y}_ {\tau}(B)
       -L^{-y}_{\hat T_1}(B)
         \right) 
          \right],
\end{eqnarray*}
 where $f_{a}(x)= x+a x^+$. Notice that ${v}_{\lambda,\beta}$ is a non negative function. Moreover, we have
\begin{lemma}\label{v_alpha}
Define
$$y_{\lambda,\beta}=-\inf\{x\in \R\;|\;v_{\lambda,\beta}(x)>0 \}. $$
We have $0<y_{\lambda,\beta}<1
                 +\lambda \beta(2+e^{\lambda})$ and   
                 \[
 \forall y< -y_{\lambda,\beta},\quad
                   v_{\lambda,\beta}(y)=0.
 \] 
\end{lemma}
We finish this paragraph with an  inequality, which will be used to derive a lower bound for the second derivative of $P$ (see the proof of the upper bound in Theorem \ref{thjd2}). \\
We define the function $C$ on $\R$ by $C(x) = x-{\lambda\beta}\ \E(B_1-x)^+$ 
and                   
we have the following  lemma,  
\begin{lemma} \label{conjecture}
For all  $x>y_{\lambda,\beta}$, we have  
$$C(x)>0.$$
\end{lemma}
 These results will be proved in \emph{ Appendix 2: A study of ${v}_{\lambda,\beta}$}.
\subsection{American put price expansion}
Throughout this section, we assume $\bar d<0$, so that $b(T)<K$.  We then have enough regularity of the American put price  to derive an expansion of $P$ around $b(T)$ along a certain parabolic branch.
\begin{theorem} \label{thjd1}
 Let $a$ be a negative  number ($a<0$) and $b(T)$ denote the limit of $b(t)$ when $t$ goes  to $T$, $b(T)=\lim_{t\rightarrow T}b(t)$. If  $\bar d<0$,  we have
  \begin{eqnarray*}
&&P(T-\theta,b(T)e^{a\sqrt{\theta}})=(K-b(T)e^{a\sqrt{\theta}})^++ C\theta^\frac{3}{2}\vv_{\lambda,\beta}(\frac{a}{\sigma})+o(\theta^\frac{3}{2}),\\
\end{eqnarray*}
where  $C=\sigma b(T) \bar{\delta} e^\lambda$, with $\lambda={\nu\{\ln {K\over b(T)}\}}$, $\bar{\delta}=\delta+\int_{y>\ln(K/b(T))}e^y\nu(dy)$  and $ \vv_{\lambda,\beta}(y)$   as defined in the previous section with  $\beta= {K\over  b(T)\bar\delta}$.
\end{theorem}
\begin{remark}\label{remarque} Notice that if $\nu$ does not charge  $\left \{\ln(\frac{K}{b(T)})\right \}$, meaning that ${{\lambda}}=0$ and $\hat T_1=\infty\;\mbox{a.s}$, then, $$\vv_{\lambda, \beta}(a)=\vv_{0}(a)=\sup_{\tau\in \mathcal{T}_{0,1}}\E\left(\int_0^\tau \left(a+ {B_{ s}}\right)ds\right).$$
In this case, the American put  price will have the  same expansion as in  the $\BS $  model, (see \cite{villeneuve1}).
\end{remark}
Before proving Theorem \ref{thjd1}, we state an elementary estimate for  
 the expectation of the  local time of  Brownian motion.
\begin{lemma}\label{lem-locB}
For all real number $a$ and for all $t>0$, we have
 \[
 0\leq \E(a-B_t)_+-a_+\leq \sqrt{t}\frac{e^{-\frac{a^2}{2t}}}{\sqrt{2\pi}}.
 \]
 \end{lemma}
 \begin{dema}\textbf{ of lemma \ref{lem-locB}:\\}{
 The first inequality follows from Jensen's inequality. For the other inequality, 
 we have
 \begin{eqnarray*}
 \E(a-B_t)_+&=&\int_{-\infty}^{a/\sqrt{t}}(a-\sqrt{t}y)e^{-y^2/2}\frac{dy}{\sqrt{2\pi}}\\
     &=&a\int_{-\infty}^{a/\sqrt{t}}e^{-y^2/2}\frac{dy}{\sqrt{2\pi}}
          +\sqrt{t}\frac{e^{-\frac{a^2}{2t}}}{\sqrt{2\pi}}
\end{eqnarray*}
Then, if $a\leq 0$, 
$$
\E(a-B_t)_+\leq \sqrt{t}\frac{e^{-\frac{a^2}{2t}}}{\sqrt{2\pi}}.
$$

If $a\geq 0$, we can write 
\begin{eqnarray*}
 \E(a-B_t)_+-a&=&-\int_{a/\sqrt{t}}^{+\infty}e^{-y^2/2}\frac{dy}{\sqrt{2\pi}}+\sqrt{t}\frac{e^{-\frac{a^2}{2t}}}{\sqrt{2\pi}}
    \leq \sqrt{t}\frac{e^{-\frac{a^2}{2t}}}{\sqrt{2\pi}}.
 \end{eqnarray*}
  }\end{dema}
  
In order to derive the expansion of the American put price, we start from the Meyer-Ito  formula (see \cite{protter}):
\begin{equation}\label{*}
(K-S_t)_+=(K-S_0)_+ +\int_0^t(-1_{\{S_s\leq K\}})
     S_s(\gamma_0 ds+\sigma dB_s)
   +\sum_{0<s\leq t}(K-S_s)_+ - (K-S_{s^-})_++{1\over 2}L^K_t,
\end{equation}
where $L^K_t$ is the local time of the process at $K$ until the date $t$.
We give, in the following lemma, an estimation of $\E L_t^K$, for small times $t$, which will allow us to neglect a part of the contribution of the  local time in the expansion of $P(t,x)$, near maturity.
\begin{lemma}\label{dam} 
Let $a $ be a negative number, $a<0$ and $ S_0=b(T)e^{a \sqrt \theta}$
. If $b(T)<K$, then we have, for all $\mathbb F$-stopping time $\tau$ with values in $[0,\theta]$,
 $$
 \E\left (L^K_\tau\right ) =2K\E\left [\left (({-a\sqrt\theta- \sigma B_{\tau}})^+- (-a\sqrt\theta-  \sigma B_{\hat T_1})^+\right )1_{\{\hat T_1<\tau\}}\right ]+o(\theta^\frac{3}{2})\leq w_0\theta^{3/2},
 $$
 where $\hat T_1 =\inf\{s\geq0\;;\;  \Delta X_s=\ln({K\over b(T)}) \}$ and $w_0$ a non-negative constant independent of $a$.
\end{lemma}
\begin{dema}\textbf{ of Lemma \ref{dam}:\\}
Let $T_1$ be the first jump time of the process $Z$ and $\tau$ an  $\mathbb F$-stopping time with values in $[0,\theta]$. We have,  the local time being a nondecreasing process,
\begin{eqnarray*}
L^K_\tau=L^K_{\tau\wedge T_1}+L^K_\tau-L^K_{\tau\wedge T_1}&=&
 L^K_{\tau\wedge T_1}+1_{\{T_1<\tau\}}\left(L^K_\tau-L^K_{T_1}\right)\\
      &\leq &L^K_{\theta\wedge T_1}+1_{\{T_1<\theta\}}\left(L^K_{T_1+\theta}-L^K_{T_1}\right).
\end{eqnarray*}

\noindent { \bf Estimating  $\E L^K_{\theta\wedge T_1}$}

In the stochastic interval $[0,T_1[$, the process $(S_t)$ matches with the process $(\check S_t)$ defined by 
$$
\check S_t=S_0e^{(\gamma_0 -\frac{\sigma^2}{2})t+\sigma B_t}.
$$
We deduce (when observing that the process $L^K$ is continuous) that
$$
L^K_{\theta\wedge T_1}=\check L^K_{\theta\wedge T_1}\leq \check L^K_\theta,
$$
where $\check L^K$ is the local time at $K$ of the process $\check S$.
Note that
$$
{1\over 2}\check L^K_\theta=(K-\check S_\theta)_+-(K-S_0)_+ -\int_0^\theta(-1_{\{\check S_s\leq K\}})
     \check S_s(\gamma_0 ds+\sigma dB_s).
$$
As the process $(\check L^K_\theta)$ increases only on $\{\check S_t=K\}$, we have 
$$
\check L^K_\theta=\check L^K_\theta1_{\{\tau_K<\theta\}},
$$
where $\check \tau_K=\inf\{t\geq 0; \check S_t> K\}$. By H\"older,
$$
\E \check L^K_\theta\leq \left(\PP(\check \tau_K<\theta)\right)^{1-\frac{1}{p}}||\check L^K_\theta||_p,\quad p>1.
$$
We easily deduce that $\E \check L^K_\theta=o(\theta^n)$, for all $n>0$.

\medskip

\noindent {\bf Estimating  $\E \left[1_{\{T_1<\tau\}}\left(L^K_{\tau}-L^K_{T_1}\right)\right]$ 
}

\noindent Notice that we have $$\E \left[1_{\{T_1<\tau\}}\left(L^K_{\tau}-L^K_{T_1}\right)\right]\leq\E \left[1_{\{T_1<\theta\}}\left(L^K_{T_1+\theta}-L^K_{T_1}\right)\right],$$
and by the strong Markov property, we obtain 
\begin{equation}\label{**}
\E \left[1_{\{T_1<\theta\}}\left(L^K_{T_1+\theta}-L^K_{T_1}\right)\right]=\E\left(1_{\{T_1<\theta\}}\E_{S_{T_1}}\!(L^K_\theta)\right),
\end{equation}
where $\E_x$ is  the expectation associated to $\PP_x$ and  $\PP_x$ defines the law of $S_t$ when $S_0=x$.
\medskip

\noindent {\bf Estimating  $\E_x \left(L^K_\theta\right)$}

\noindent Let $T_1$ be the first jump time of the process $Z$.
We then have
\begin{eqnarray*}
L^K_\theta=L^K_{\theta\wedge T_1}+L^K_\theta-L^K_{\theta\wedge T_1}
\end{eqnarray*}
According to equality (\ref{*}) we  deduce, using the compensation formula (see \cite{bertoin})
\begin{eqnarray*}
{1\over 2}L^K_\theta&=&
(K-S_\theta)_+-(K-S_0)_+ +\int_0^\theta1_{\{S_{s^-}\leq K\}}
     S_{s^-}(\gamma_0 ds+\sigma dB_s)\\
   &&-\int_0^\theta ds\int \Phi(S_{s^-},y) \nu(dy)+M_\theta,
\end{eqnarray*}
where  $\Phi(x,y)=(K-xe^y)_+-(K-x)_+$ and $(M_t)$ is a martingale which vanishes at $0$.
Taking expectations, we have
     $$
     {1\over 2}\E\left(L^K_{\theta}\right)
     =\E(K-S_{\theta})_+-(K-S_0)_+ +\E\int_0^{\theta}
        \left(\gamma_0 S_s1_{\{S_s\leq K\}}
   -\int \Phi(S_s,y)  \nu(dy)\right)ds.
     $$
We deduce easily from this equality that 
$$
{1\over 2}\E_x\left(L^K_\theta\right)=\E_x(K-S_\theta)_+-(K-x)_+ +xO(\theta)
$$
with $O(\theta)$ independent of  $x$.
We have 
\begin{eqnarray*}
\E_x(K-S_\theta)_+-(K-x)_+ & = & \E_x(K-xe^{(r-\delta-\frac{\sigma^2}{2})\theta+\sigma B_\theta+\tilde{Z}_\theta})_+-(K-x)_+  
\end{eqnarray*}
We also have 
\begin{eqnarray*}
 \E \left|e^{(r-\delta-\frac{\sigma^2}{2})\theta+\sigma B_\theta+\tilde{Z}_\theta}  -e^{\sigma B_\theta}\right|
         &=&e^{\sigma^2\theta/2}
            \E \left|e^{(r-\delta-\frac{\sigma^2}{2})\theta+\tilde{Z}_\theta}  -1\right|\\
            &=&O(\theta)
\end{eqnarray*}
Therefore
\begin{eqnarray*}
\E_x(K-S_\theta)_+-(K-x)_+ & = & \E(K-xe^{\sigma B_\theta})_+-(K-x)_+  +xO(\theta)\\
            &=&\E(K-x(1+\sigma B_\theta))_+-(K-x)_+  +xO(\theta)\\
            &=&x\sigma\left(\E\left(\frac{K-x}{x\sigma}- B_\theta\right)_+-\left(\frac{K-x}{x\sigma}\right)_+\right)+
               xO(\theta).
\end{eqnarray*}
Hence, using lemma ~\ref{lem-locB} above,
$$
\E_x(K-S_\theta)_+-(K-x)_+ \leq
       x\sigma\sqrt{\theta/(2\pi)}\exp\left(-\frac{(K-x)^2}{2x^2\sigma^2\theta}\right)+xO(\theta).
$$
Going back  to (\ref{**}), we obtain
\begin{eqnarray*}
\lefteqn{{1\over 2} \E \left[1_{\{T_1<\theta\}}\left(L^K_{T_1+\theta}-L^K_{T_1}\right)\right]}\\&\leq & 
    \sigma\sqrt{\frac{\theta}{2\pi}}\E\left(1_{\{T_1<\theta\}} S_{T_1}
                 \exp\left(-\frac{(K-S_{T_1})^2}{2S_{T_1}^2\sigma^2 \theta}\right)\right)+\E\left(1_{\{T_1<\theta\}} S_{T_1}\right)O(\theta)
                 \\
                 &=&\sigma\sqrt{\frac{\theta}{2\pi}}S_0
                 \E\left(1_{\{T_1<\theta\}} e^{(\gamma_0-\frac{\sigma^2}{2})T_1+\sigma B_{T_1}+Z_{T_1}}
                 \exp\left(-\frac{(K-S_{T_1})^2}{2S_{T_1}^2\sigma^2 \theta}\right)\right)+O(\theta^2).
\end{eqnarray*}
At this stage, we notice that $\PP(T_1\leq \theta)=1-e^{-\lambda \theta}=O(\theta)$ and that, conditionally on $\{T_1\leq \theta\}$, $T_1$ is uniformly distributed  on $[0,\theta]$.

 As $Z_{T_1}$ 
is independent of both $T_1$ and $B$, we see that, conditionally to $\{T_1<\theta\}$,
${S_{T_1}}$ has the same law as 
$$K\exp\left \{\left ({V-\ln( \frac{K}{b(T)} )}\right )+\sqrt \theta\left (a+(\gamma_0-\frac{\sigma^2}{2}) \sqrt \theta U+\sigma g \sqrt{U}\right )\right \},$$
where $U$, $g$ and $V$ are two independent random variables, $U$ is uniform on $[0,1]$,
$g$ standard Gaussian  and $V$ has the same law as $Z_{T_1}$.
Therefore, we can state that there exists a  non negative constant independent of $a$ such that
$$ \E(L^K_\theta\;|\; S_0=b(T)e^{a \sqrt \theta})\leq w_0 \theta^{3/2}=O(\theta^{3/2}).$$

\noindent {\bf Estimating $\E \left[L^K_\tau-L^K_{\tau\wedge T_1}\right]$, in the case where $\nu\{\ln({K\over b(T)})\}=0$:}

 If we assume $\nu\{\ln({K\over b(T)})\}=0$, which means that $V-\ln( \frac{K}{b(T)})\neq 0$ a.s , we obtain, by dominated convergence that 
$$
\lim_{t\downarrow 0}\E\left(S_{T_1}
                 \exp\left(-\frac{K^2(1-{S_{T_1}\over K})^2}{2S_{T_1}^2\sigma^2 \theta}\right)\;|\; T_1<\theta\right)=0.
$$
Therefore $\E \left[1_{\{T_1<\theta\}}\left(L^K_{T_1+\theta}-L^K_{T_1}\right)\right]=o(\theta^{3/2})$, 
hence  $$\E \left(L^K_{\theta}\;|\; S_0=b(T)e^{a \sqrt \theta}\right)=o(\theta^{3/2})$$

\noindent {\bf Estimating	 $\E \left[L^K_\tau-L^K_{\tau\wedge T_1}\right]$, in the case where $\nu\{\ln({K\over b(T)})\}>0$:}

Let us introduce the processes $\hat X$ and $\hat Z$ such that
$$\hat Z_t=\sum_{s<t}\Delta \tilde X_s 1_{\{\Delta\tilde  X_s=\ln\frac{K}{b(T)}\}}\quad\mbox{and}\quad \hat X=\tilde X-\hat Z,$$ and $\hat T_1 =\inf\{s\geq0,\quad \hat Z_t\neq 0 \}$.
Then, since $\tau\leq \theta$, we have $$\E \left[L^K_{\tau\wedge \hat T_1}-L^K_{\tau\wedge T_1}\right]=\E \left[(L^K_{\tau\wedge \hat T_1}-L^K_{\tau\wedge T_1})1_{\{T_1<\tau\wedge \hat T_1\}}\right]=o(\theta^\frac{3}{2}).$$ 
Indeed,  on $\{\tau<\hat T_1\}$, the process $\tilde X$ matches with the process $\hat X$ whose Lévy measure does not charge the point $\{\ln({K\over b(T)})\}$, (we are in the same case as $\nu\{\ln({K\over b(T)})\}=0$). And on $\{T_1<\hat T_1\leq \tau\}\subset\{T_1<\hat T_1\leq \theta\}$,  the process $Z$ has jumped two times   before $\theta$, however,  $\PP\left (\sum_{s\leq \theta}1_{\{\Delta Z_s\neq 0\}}\geq 2\right )=O(\theta^2)$. Thus,
$$\E \left[L^K_\tau-L^K_{\tau\wedge T_1}\right]=\E \left[L^K_\tau-L^K_{\tau\wedge \hat T_1}\right]+o(\theta^{3/2}).$$
Besides,
\begin{eqnarray*}
\lefteqn{{1\over 2}\E \left[L^K_\tau-L^K_{\tau\wedge \hat T_1}\right]={1\over 2}\E \left[(L^K_\tau-L^K_{\tau\wedge \hat T_1})1_{\hat T_1<\tau}\right]}\\
&= &\E\left [\left ((K-S_\tau)^+- (K-S_{ \hat T_1})^+\right )1_{\hat T_1<\tau}\right ]+o(\theta^{3/2})\\
&= & \E\left [\left ((K-{S_0 K\over b(T)}e^{\hat X_{\hat T_1}-\tilde X_{\hat T_1}+\tilde X_\tau})^+- (K-{S_0 K\over b(T)}e^{\hat X_{\hat T_1}})^+\right )1_{\hat T_1<\tau}\right ]+o(\theta^\frac{3}{2}).
\end{eqnarray*}
Since $\PP\left (\sum_{s\leq \theta}1_{\{\Delta \tilde X _s\neq 0\}}\geq 2\right )=O(\theta^2)$,  conditionally on $\{\hat T_1<\tau\}$, we can assume  that $N_\theta=1$, where $N_\theta$ denotes the number of jumps of $\tilde X$ up to $\theta$, $N_\theta =\sum_{s\leq \theta}1_{\{\Delta \tilde X _s\neq 0\}}=\sum_{s\leq \theta}1_{\{\Delta Z _s\neq 0\}}$. Noticing that ${S_0 K\over b(T)}= Ke^{a \sqrt\theta}$, we obtain
\begin{eqnarray}
\lefteqn{\E \left[L^K_\tau-L^K_{\tau\wedge \hat T_1}\right]}\nonumber\\&=&2\E\left [\left ((K-{S_0 K\over b(T)}e^{\hat X_{\hat T_1}-\tilde X_{\hat T_1}+\tilde X_\tau})^+- (K-{S_0 K\over b(T)}e^{\hat X_{\hat T_1}})^+\right )1_{\{N_\theta=1\}}1_{\{\hat T_1<\tau\}}\right ]+o(\theta^\frac{3}{2})\nonumber\\
&= &2K\E\left [\left (({-a\sqrt\theta- \mu\tau- \sigma B_{\tau}})^+- (-a\sqrt\theta- \mu{\hat T_1}- \sigma B_{\hat T_1})^+\right )1_{\{\hat T_1<\tau\}}\right ]+o(\theta^\frac{3}{2})\nonumber
 \\
&= &2K\E\left [\left (({-a\sqrt\theta- \sigma B_{\tau}})^+- (-a\sqrt\theta- \sigma B_{\hat T_1})^+\right )1_{\{\hat T_1<\tau\}}\right ]+o(\theta^\frac{3}{2})\nonumber
\end{eqnarray}
The  last two equalities follow from $\PP(\hat T_1<\tau)=O(\theta)$, $\left |(1-e^x + x)1_{\{x\leq 0\}}\right |\leq \frac{x^2}{2}$ and the fact that, for all stropping time $\varrho$ with values in $[0,\theta]$, we have $$\theta\E\left (a+ \mu {\varrho\over\sqrt \theta}+ { \sigma \over \sqrt\theta } B_{\varrho}\right )^2\leq C \theta .$$

\end{dema}
\begin{dema} \textbf{of Theorem \ref{thjd1}:}
First of all, we recall our notation $\check{X}_t=\tilde{X}_t-Z_t$, $\check{S}_t=\tilde{S}_t/e^{Z_t}$ (i.e the continuous part of the processes) and $T_1$ the first jump time $T_1=\inf\{t>0|Z_t\neq0\}$ and from now on, we consider $S_0$ as a function of $\theta$. More precisely, we denote by $S_0^\theta=b(T)e^{a\sqrt\theta}=e^{x_0+a\sqrt\theta}$, with $a<0$ and $x_0=\ln (b(T))$.\\
According to equation (\ref{*}), we have for all stopping times $\tau \in \mathcal{T}_{0,\theta}$,
\begin{eqnarray}\lefteqn{\E\left [e^{-r\tau}(K-S_\tau)_+\right ]-(K-S_0)^+}\nonumber\\&=& 
\E\left[\int_0^\tau\left(e^{-rs}1_{\{S_s\leq K\}}
\left(-rK+\delta S_s+S_s\int (e^y-1)\nu(du)\right)\right.\right.\nonumber\\
&&\left.\left.+ e^{-rs}\int \left [(K-S_se^y)^+-(K-S_s)^+\right ]) \nu(dy) \right)ds\right]+{1\over 2}\E\left(\int_0^\tau e^{-rs}dL^K_s\right)\nonumber
\\&=& \I^a(\tau)+\J^a(\tau), \label{dl europ}
\end{eqnarray}
where
\begin{eqnarray*}
\I^a(\tau)&=&\E\left[\int_0^\tau\left(e^{-rs}1_{\{S_s\leq K\}}
\left(-rK+\delta S_s+S_s\int (e^y-1)\nu(du)\right)\right.\right.\\
&&\left.\left.+ e^{-rs}\int \left [(K-S_se^y)^+-(K-S_s)^+\right ]) \nu(dy) \right)ds\right]
\end{eqnarray*}
and 
$$\J^a(\tau)={1\over 2}\E\left(\int_0^\tau e^{-rs}dL^K_s\right).
$$
At this stage, since $\J^a\geq 0$, we can state that , given $S_0=b(T)e^{\lambda \sqrt \theta}$, we have 
\begin{eqnarray}
\I^a(\tau) \leq \E\left [e^{-r\tau}(K-S_\tau)_+\right ]-(K-S_0)^+= \I^a(\tau)+\J^a(\tau)\leq \I^a(\tau)+w_0\theta^{3\over 2},
\label{encadrement P_e}
\end{eqnarray}
the last inequality follows from Lemma \ref{dam}.
In what follows, we will express $\I^a$ and $\J^a$ in  more appropriate forms. Let us  start with $\J^a$. 

\medskip

\noindent {\bf Estimating  $\J^a$:}\\
Recall that $\hat T_1=\inf \{t\geq0\;;\;\Delta \tilde X_t=\ln{K\over b(T)}\}$ and set$S_0=b(T)e^{\lambda \sqrt \theta}$ with  $\lambda<0$, then according to Lemma \ref{dam}, we have $
 \E\left (L^K_\theta\right )=O(\theta^{3/2}),
 $   therefore
\begin{eqnarray}
\J^a(\tau)&=& {1\over 2}\E\left(L^K_\tau\right)+o(\theta^\frac{3}{2}) \nonumber
 \\
&= &K\E\left [\left (({-a\sqrt\theta- \sigma B_{\tau}})^+- (-a\sqrt\theta- \sigma B_{\hat T_1})^+\right )1_{\{\hat T_1<\tau\}}\right ]+o(\theta^\frac{3}{2}).\label{derna3}
\end{eqnarray}
\medskip

\noindent {\bf Estimating of $\I^a$:}\\
\medskip
First of all, remark that we have
\begin{eqnarray}
\lefteqn{\E\left[\int_0^\tau\left(e^{-rs}1_{\{S_s > K\}}
\int \left [(K-S_se^y)^+\hspace{-2MM}-(K-S_s)^+\right ] \nu(dy) \right)ds\right]\nonumber}&&\\
&\leq & K \nu(\R) \int_0^\theta \PP\{S_s > K\}ds \nonumber \\
&\leq&K \nu(\R)\int_0^\theta \PP\{S_s > K,T_1>\theta\}+\PP\{S_s > K,T_1\leq\theta\}ds\label{CondT1} \\
&\leq&K \nu(\R)\left (\int_0^\theta \PP\{\check S_s >K\}ds+\theta \PP\{T_1\leq\theta\}\right )=O(\theta^2). \nonumber
\end{eqnarray}
And noticing that $$1_{\{x\leq K\}}\left (x(e^y-1)+\left[(K-xe^y)^+-(K-x)^+\right ]\right )=(xe^y-K)^+1_{\{x\leq K\}},$$
 we thus obtain 
\begin{eqnarray*}
\I^a(\tau)=\E\left(\int_0^\tau \hspace{-3MM}e^{-rs}1_{\{S_s\leq K\}}\left(-rK+\delta S_s+\int (S_se^y-K)^+ \nu(dy)\right)ds\right)+o(\theta^\frac{3}{2}).
\end{eqnarray*}
We can also omit  $e^{-rs}$ in the expression as an error  of the order of  $O(\theta^2)$.
Then we  obtain, for all  stopping times $\tau$ with values in $[0,\theta]$ 
\begin{eqnarray*}
\I^a(\tau)=\E\left(\int_0^\tau \hspace{-3MM}1_{\{S_s\leq K\}}\left(-rK+\delta S_s+\int (S_se^y-K)^+ \nu(dy)\right)ds\right)+o(\theta^\frac{3}{2}).
\end{eqnarray*}
We  denote 
$$h(x)=-rK+\delta e^{x}+\int (e^{x}e^y-K)^+ \nu(dy),$$
and recall that $S_t=S_0^\theta e^{\tilde X_t}=b(T)e^{a\sqrt\theta+\tilde X_t}=b(T)e^{\tilde X^{a\sqrt\theta}_t}=e^{x_0+\tilde X^{a\sqrt\theta}_t},$
where $\tilde {X}^{y}_t=y+\tilde X_t$  .
We thus have
\begin{eqnarray}
\I^a(\tau)&=&\underbrace{\E\left(\int_0^\tau\hspace{-3MM}
1_{\{a\sqrt{\theta}+\tilde{X}_s\leq \ln \frac{K}{b(T)} \}} h({x_0+a\sqrt{\theta}+\tilde{X}_s})ds\right)}_{(I)} +o(\theta^\frac{3}{2}).
\label{eq2}
\end{eqnarray}
Now, we will try to  express  the quantity ($I$) under a more appropriate form.
The first step is to neglect the contribution of the finite variation part of the process $\tilde X$.
Notice that
$$
\left |1_{\{x\leq \ln(K)\}}h(x)\right |\leq K(r\vee |\bar d|)\quad\text{and }\quad |h(x)-h(y)|\leq |e^x-e^y|\left (\delta+\int_{y>0} e^y\nu(du)\right ).
$$  Moreover, for all 
$(x,y)\in \R^2$, we have 
\begin{eqnarray*}
\lefteqn{\left |1_{\{x\leq \ln(K) \}} h(x)-1_{\{ y\leq \ln(K)\}} h(y)\right |}\\
&=&\left |
\left(h(x)-h(y)\right ) 1_{\{x\vee  y\leq \ln(K) \}} +h(x)1_{\{x\leq \ln(K)< y \}}
-h(y)1_{\{ y\leq \ln(K)<x\}} 
\right |\\
&\leq&
A_0\left| e^x-e^y\right | 1_{\{x\vee  y\leq \ln(K) \}}A_1 \left (1_{\{ \ln(K)< y \}}
+1_{\{  \ln(K)<x\}} 
\right ),
\end{eqnarray*}
where $A_1=K(r\vee|\bar d|)>0$ and $A_0=\delta+\int_{y>0} e^y\nu(du)$.
Let   $k_b=\ln\left (\frac{K}{b(T)}\right )>0$ and recall that $ \tilde X_t-\sigma B_{t }=(\gamma_0-\frac{\sigma^2}{2}) t+Z_t$, then 
\begin{eqnarray*}
\lefteqn{\left |1_{\{x_0+a\sqrt{\theta}+\tilde{X}_s\leq \ln K \}} h({x_0+a\sqrt{\theta}+\tilde{X}_s})-1_{\{ x_0+{a\sqrt{\theta}}+\sigma B_{s }\leq \ln K \}} h({x_0+ {a\sqrt{\theta}}+\sigma B_{s}})\right |}
\\
&\leq &
A_0\left | e^{x_0+a\sqrt{\theta}+\tilde{X}_s}-e^{x_0+{a\sqrt{\theta}}+\sigma B_{s }}\right | 
1_{\{
\tilde{X}_s \vee \sigma B_{s} \leq k_b-a\sqrt{\theta} 
\}} 
+C(1_{\{ k_b-{a\sqrt{\theta}}< \sigma B_{s} \}}+1_{\{ k_b-{a\sqrt{\theta}}<\tilde{X}_s \}})
\\
&\leq &
A_0b(T)e^{\sigma B_s}\left | e^{(\gamma_0-\frac{\sigma^2}{2}) s+Z_s}-1 \right | 
+
C(1_{\{ k_b< \sigma B_{s} \}}+1_{\{ k_b<\tilde{X}_s \}}),
\end{eqnarray*}
the last inequality is due to $a<0$ and $e^x_0=B(T)$.\\
Taking the expectation, we obtain, for all $ s\in[0, \theta]$
\begin{eqnarray*}
\PP(k_b< \sigma B_{s})\leq \PP(\frac{k_b}{\sigma\sqrt{\theta}}< \sigma B_1)\leq C \sqrt{\theta}e^{-\frac{k_b^2}{2\sigma^2\theta}},
\end{eqnarray*}
for  $\theta$ small enough, we  have $\frac{k_b}{2}<k_b-(\gamma_0-\frac{\sigma^2}{2}) s$,  then
\begin{eqnarray*}
\PP(k_b< \tilde X_{s})\leq \PP(\frac{k_b-(\gamma_0-\frac{\sigma^2}{2}) s}{\sigma\sqrt{\theta}}<  B_1)+\PP(T_1\leq\theta)\leq C \sqrt{\theta}e^{-\frac{k_b^2}{8\sigma^2\theta}}+A\theta
\end{eqnarray*}
and \[\E\left (e^{\sigma B_s}\left | e^{(\gamma_0-\frac{\sigma^2}{2}) s+Z_s}-1 \right | \right )\leq e^{\frac{\sigma^2}{2} s}\left | e^{(\gamma_0-\frac{\sigma^2}{2}) s}-1 \right |+ e^{\gamma_0 s}\E\left | e^{Z_s}-1 \right | \leq D\theta.\]
Hence, 
\begin{eqnarray*}
\lefteqn{\int_0^\theta\E \left(\left |1_{\{x_0+a\sqrt{\theta}+\tilde{X}_s\leq \ln K \}} h({x_0+a\sqrt{\theta}+\tilde{X}_s})\right .\right .}
\\&&-\left .\left .1_{\{ x_0+{a\sqrt{\theta}}+\sigma B_{s }\leq \ln K \}} h({x_0+ {a\sqrt{\theta}}+\sigma B_{s}})\right |\right )ds=O(\theta^2).
\end{eqnarray*}
Thanks to this estimation,  equation (\ref{eq2})  becomes 
\begin{eqnarray}
\I^a(\tau)&=& \E\left(\int_0^\tau\hspace{-3MM}
1_{\{{a\sqrt{\theta}}+\sigma B_s\leq \ln \frac{K}{b(T)} \}} h({x_0+{a\sqrt{\theta}}+\sigma B_s})ds\right) +o(\theta^\frac{3}{2}).\label{eq3}
\end{eqnarray}
The function $h$ is convex, therefore it is right and left differentiable. Particularly, we have all $x<\ln(K)$,
$$h'_g(x)=e^x\left (\delta +\int e^{y}1_{\{y>\ln(K)-x\}}\nu(dy)\right )$$ and $$h'_d(x)=e^x\left (\delta +\int e^{y}1_{\{y \geq \ln(K)-x\}}\nu(dy)\right ).$$ 
Hence, we can write 
\begin{eqnarray*}
h'_d(x_0)(x-x_0)^+-h'_g(x_0)(x-x_0)^- \leq h(x)-h(x_0)\leq h'_g(x)(x-x_0)^+-h'_d(x)(x-x_0)^-,
\end{eqnarray*}
hence
\begin{eqnarray*}
0&\leq& h(x)-\left (h(x_0)+h'_d(x_0)(x-x_0)^+-h'_g(x_0)(x-x_0)^-\right )
\\&\leq&\left ( h'_g(x)-h'_d(x_0)\right )(x-x_0)^++\left (h'_g(x_0) -h'_d(x)\right )(x-x_0)^-
\\&=&\left ( h'_g(x\vee x_0)-h'_d(x\wedge x_0)\right )|x-x_0|.
\end{eqnarray*}
Thanks to the equation characterizing $b(T)$ when $\bar d<0$, we have  $h(x_0)=h(\ln(b(T))=0$. We thus obtain, by setting  $\Delta h'(x_0)=h'_d(x_0)- h'_g(x_0)$,
\begin{eqnarray*}
 h(x_0+x)
 =\Delta h'(x_0) x^+ + h'_g(x_0)x+\left| x\right|\tilde{R}( x), 
\end{eqnarray*} 
where $ \tilde{R}(x)\longrightarrow_{x\rightarrow 0}0$, and 
\begin{eqnarray*}
0\leq \tilde{R}(x) 
 &\leq &\left (h'_g(x_0+x^+)-h'_d(x_0-x^-)
\right)
\\
&\leq& L \left (1+e^{x}\right ),
\end{eqnarray*} 
 with $L$ a positive constant.
We can then write
\begin{eqnarray}
 \lefteqn{1_{\{a\sqrt{\theta}+\sigma B_s\leq \ln \frac{K}{b(T)}\}} h(x_0+a\sqrt{\theta}+\sigma B_s)}\nonumber&&\\&=&\left (\Delta h'(x_0)(a\sqrt{\theta}+\sigma B_s)^+ +h'_g(x_0)a\sqrt{\theta}+\sigma B_s\right )
 \left (1-1_{\{a\sqrt{\theta}+\sigma B_s> \ln \frac{K}{b(T)} \}}\right )\nonumber\\&&+\left| a\sqrt{\theta}+\sigma B_s\right|\tilde{R}( a\sqrt{\theta}+\sigma B_s)1_{\{\tilde{X}^{a
 \sqrt{\theta}}_s\leq \ln \frac{K}{b(T)} \}}. \label{1_h}
\end{eqnarray}  
We state that
\begin{eqnarray}
\left |\E \int_0^\tau \left| a\sqrt{\theta}+\sigma B_s\right|\tilde{R}( a\sqrt{\theta}+\sigma B_s)1_{\{\tilde{X}^{a
 \sqrt{\theta}}_s\leq \ln \frac{K}{b(T)} \}} ds\right |=o(\theta^\frac{3}{2})\label{R1}\\
\left |\E \int_0^\tau \left (\Delta h'(x_0)(a\sqrt{\theta}+\sigma B_s)^+ +h'_g(x_0)a\sqrt{\theta}+\sigma B_s\right )1_{\{a\sqrt{\theta}+\sigma B_s> \ln \frac{K}{b(T)} \}}ds\right |=o(\theta^\frac{3}{2}),\nonumber\\\label{R2}
\end{eqnarray}
Indeed, we have  for (\ref{R1}),  by setting $s=u \theta$,
\begin{eqnarray*}
\lefteqn{\left |\E\left(\int_0^\tau \left| a\sqrt{\theta}+\sigma B_s\right|\tilde{R}( a\sqrt{\theta}+\sigma B_s)1_{\{a\sqrt{\theta}+\sigma B_s\leq \ln \frac{K}{b(T)} \}}ds\right)\right |}&&\\
&= &\theta^\frac{3}{2}\int_0^1 \E\left[\left|
a+\sigma B_{ s}\right| \tilde{R}( \sqrt\theta(a+\sigma B_{ s}))1_{\{a\sqrt\theta+\sigma\sqrt\theta B_{ s}\leq \ln \frac{K}{b(T)} \}}\right]ds .
\end{eqnarray*}
As $|\tilde{R}(x)|\leq L(e^x+1)$ and $|\tilde{R}(x)|\longrightarrow_{x \rightarrow 0}  0$, we have by bounded convergence
\begin{eqnarray}
\int_0^1 \E\left[\left
|a+\sigma B_{ s}\right| 
\tilde{R}( \tilde{X}_{\theta s}^{a\sqrt{\theta}})1_{\{a
 \sqrt{\theta}+\sigma B_s\leq \ln \frac{K}{b(T)} \}}\right]ds
&\longrightarrow_{\theta\rightarrow 0} & 0 .\nonumber
\end{eqnarray}
And for the estimate in (\ref{R2}), we have
\begin{eqnarray*}
\lefteqn{\left |\E \int_0^\tau \left (\Delta h'(x_0)(a\sqrt{\theta}+\sigma B_s)^+ +h'_g(x_0)a\sqrt{\theta}+\sigma B_s\right )1_{\{a\sqrt{\theta}+\sigma B_s> \ln \frac{K}{b(T)} \}}ds\right |}&&\\
&\leq & C\sqrt{\theta}\int_0^\theta \E \left [\left (\left| a\right|+\sigma\sqrt{\frac{s}{\theta}}\left | B_1 \right |\right ) 1_{\{a+\sigma B_1>\frac{1}{\sqrt{\theta}} \ln \frac{K}{b(T)}\}} \right ]ds\\
&\leq & C\theta^\frac{3}{2}\sqrt{ \E \left (\left| a\right |+\left |B_1 \right |\right  )^2}\sqrt{ \PP{\{a+\sigma B_1>\frac{1}{\sqrt{\theta}} \ln \frac{K}{b(T)}\}}}
\\
&=&O(\theta^n).
\end{eqnarray*}
Therefore, taking  the expectation of the integral of \eqref{1_h} between $0$ and all stopping time $\tau \in \mathcal{T}_{0,\theta}$  gives 
\begin{eqnarray}
\lefteqn{\I^a(\tau)=\E\int_0^\tau 1_{\{a\sqrt{\theta}+\sigma B_s\leq \ln \frac{K}{b(T)}\}}h(x_0+a\sqrt{\theta}+\sigma B_s )ds\nonumber}&&\\
&=&
 h'_g(x_0)\E\int_0 ^\tau \hspace{-2mm} \left (a\sqrt{\theta}+\sigma{B}_s \right )ds +\Delta h'(x_0)\E\int_0^\tau \hspace{-2mm} \left (a\sqrt{\theta}+\sigma{B}_s\right )^+ \hspace{-2mm}ds + o(\theta^\frac{3}{2}),\nonumber\\
 &=&
b(T)\bar \delta \E\int_0 ^\tau \hspace{-2mm} \left (a\sqrt{\theta}+\sigma{B}_s \right )+{\lambda\beta} \left (a\sqrt{\theta}+\sigma{B}_s\right )^+ \hspace{-2mm}ds + o(\theta^\frac{3}{2}),\label{derna2}
 \end{eqnarray}
 with $\bar\delta=\delta+\int_{y>\ln{K\over b(T)}}\hspace{0mm} e^y\nu(dy)$ , ${\beta}={K \over b(T) \bar\delta}$, ${\lambda}=\nu\{\ln{K\over b(T)}\}$ and we recall that
 $h'_g(x_0)=b(T) \bar\delta$ {and} $\Delta h'(x_0)=K\nu\{\ln{K\over b(T)}\}$ then ($\lambda \beta={\Delta h'(x_0)\over h'_g(x_0)}$).
\\Comming back to   \eqref{dl europ} and using \eqref{derna3} and \eqref{derna2}, we obtain
\begin{eqnarray*}
\E\left(e^{-r\tau} (K-S_\tau)^+\right)&=&(K-S_0)^+
      +     \E\left(b(T)\bar\delta \int_0^{\tau} \left(
   a\sqrt{\theta}+ \sigma B_s
     +{\lambda\beta}(a\sqrt{\theta}+ \sigma B_s)^+ 
     \right)ds\right)\\
     &&\left. +
K  1_{\{\hat T_1<\tau\}}
       \left((a\sqrt{\theta}+\sigma B_\tau)^+
       -(a\sqrt{\theta}+\sigma B_{\hat T_1})^+
         \right)
         \right)
         +o(\theta^{3/2}),
\end{eqnarray*}
with  $o(\theta^{3/2})$ independent of  $\tau$. Hence
\begin{eqnarray*}
P(T-\theta,b(T)e^{a\sqrt{\theta}})&=&
   (K-b(T)e^{a\sqrt{\theta}})^+
   +\sigma b(T)\bar\delta
  \bar v_{\lambda, \beta, \theta}(a/\sigma)+o(\theta^{3/2}), 
   \end{eqnarray*}
 where   $\bar{v}_{\lambda,\beta,\theta}$
   defined by
\[
\bar{v}_{\lambda,\beta, \theta}(y)=\sup_{\tau\in \T_{0,\theta}}
\E\left(\int_0^{\tau} f_{\lambda\beta}(
   y\sqrt{\theta}+ B_s
    )ds
     +\beta 
      1_{\{\hat T_1<\tau\}}
       \left((y\sqrt{\theta}+ B_\tau)^+
       -(y\sqrt{\theta}+ B_{\hat T_1})^+
         \right)\right),
\]
with $f_{a}(x)=x+a x^+$. To simplify the expression of 
$\bar{v}_{\lambda,\beta, \theta}$, we notice first that, if we set  $B^\theta_t=B_{\theta t}/\sqrt{\theta}$, we can  write
\begin{eqnarray*}
\bar{v}_{\lambda,\beta, \theta}&=&
          \sqrt{\theta}
          \sup_{\tau\in \T_{0,\theta}}
\E\left(\int_0^{\tau} f_{\lambda\beta}(
   y+ B^\theta_{s/\theta}
    )ds
     +\beta 
      1_{\{\hat T_1<\tau\}}
       \left((y+ B^\theta_{\tau/\theta})^+
       -(y+ B^\theta_{\hat T_1/\theta})^+
         \right)\right)\\
         &=&
  \sqrt{\theta}
          \sup_{\tau\in \T_{0,\theta}}
\E\left(\theta\int_0^{\tau/\theta} f_{\lambda\beta}(
   y+ B^\theta_{s}
    )ds
     +\beta 
      1_{\{\hat T_1<\tau\}}
       \left((y+ B^\theta_{\tau/\theta})^+
       -(y+ B^\theta_{\hat T_1/\theta})^+
         \right)\right)       
\end{eqnarray*}
We also notice that $\tau\in \T_{0,\theta}$ if and only if  $\tau/\theta\in \T_{0,1}^\theta$, where $\T_{0,1}^\theta$ is the set of the stopping times of the filtration $(\F_{\theta t})_{t\geq 0}$, with values in $[0,1]$, then
\[
\bar{v}_{\lambda,\beta, \theta}=
\sqrt{\theta}
          \sup_{\tau\in \T_{0,1}^\theta}
\E\left(\theta\int_0^{\tau} f_{\lambda\beta}(
   y+ B^\theta_{s}
    )ds
     +\beta 
      1_{\{\hat T_1<\theta\tau\}}
       \left((y+ B^\theta_{\tau})^+
       -(y+ B^\theta_{\hat T_1/\theta})^+
         \right)\right) 
\]
 Note that 
$\bar{v}_{\lambda,\beta, \theta}(y)$ does not  change  if we replace
$\T_{0,1}^\theta$ by $\hat\T_{0,1}$ the set of the stopping times of the natural filtration of the  couple $(B^\theta_t, \hat N_{\theta t})$, where 
$\hat N$ is defined by  
\[
\hat N_t=\sum_{0<s\leq t} 1_{\{\Delta Z_s=\ln(K/b(T))\}}.
\]
The  processes $(\hat N_{\theta t})_{t\geq 0}$ is a Poisson process
with intensity $\theta \lambda$, where $\lambda=\nu\{\ln( K/b(T))\}$.
Under the  probability $\hat \PP$, defined by
\[
\frac{d\hat \PP}{d\PP}=\theta^{\hat N_1}e^{-\lambda(\theta-1)},
\] 
the process $(B_t, \hat N_t)_{0\leq t\leq 1}$
  has the same law as $(B^\theta_t, \hat N_{\theta t})_{0\leq t\leq 1}$.
Hence,
\begin{eqnarray*}
\bar{v}_{\lambda,\beta, \theta}(y)&=&
\sqrt{\theta}
          \sup_{\tau\in \T_{0,1}}
\E\left[\theta^{\hat N_1}e^{-\lambda(\theta-1)}
\left(\theta\int_0^{\tau} f_{\lambda\beta}(
   y+ B_{s}
    )ds
     +\beta 
      1_{\{\hat T_1<\tau\}}
       \left((y+ B_{\tau})^+
       -(y+ B_{\hat T_1})^+
         \right)\right) \right]\\
         &=&
        \sqrt{\theta}
          \sup_{\tau\in \T_{0,1}}
\E\left[\theta^{\hat N_\tau}e^{-\lambda\tau (\theta-1)}
\left(\theta\int_0^{\tau} f_{\lambda\beta}(
   y+ B_{s}
    )ds
     +\frac{\beta }{2}
      1_{\{\hat T_1<\tau\}}
       \left(L^{-y}_ {\tau}(B)
       -L^{-y}_{\hat T_1}(B)
         \right)\right) \right], 
\end{eqnarray*}
where $L^{-y}(B)$ denotes the local time of  $B$ at $-y$.
We have for $\tau\in \T_{0,1}$,
\begin{eqnarray*}
\E\left[\theta^{\hat N_\tau}e^{\lambda\tau (\theta-1)}
\left(\theta\int_0^{\tau} f_{\lambda\beta}(
   y+ B_{s}
    )ds\right)
\right]=\theta\E\left[ 1_{\{\hat N_\tau=0\}}
        e^{-\lambda\tau (\theta-1)}
\left(\int_0^{\tau} f_{\lambda\beta}(
   y+ B_{s}
    )ds\right)
\right]+\theta R_\tau,
\end{eqnarray*}
and if $\theta\leq 1$
\begin{eqnarray*}
|R_\tau|&\leq &\theta\E\left[ 1_{\{\hat N_\tau\geq 1\}}
       e^{-\lambda\tau (\theta-1)}
\left(\int_0^{1} |f_{\lambda\beta}(
   y+ B_{s}
    )|ds\right)
\right]=O(\theta).
\end{eqnarray*}
Hence,
\begin{eqnarray*}
\E\left[\theta^{\hat N_\tau}e^{-\lambda\tau (\theta-1)}
\left(\theta\int_0^{\tau} f_{\lambda\beta}(
   y+ B_{s}
    )ds\right)
\right]=\theta\E\left[ 1_{\{\hat N_\tau=0\}}
        e^{\lambda\tau }
\left(\int_0^{\tau} f_{\lambda\beta}(
   y+ B_{s}
    )ds\right)
\right]+O(\theta^2),
\end{eqnarray*}
Besides,
\begin{eqnarray*}
\E\left[\theta^{\hat N_1}e^{-\lambda(\theta-1)}
      1_{\{\hat T_1<\tau\}}
       \left((L^{-y}_ {\tau}(B)
       -L^{-y}_{\hat T_1}(B)
                \right) \right]&=&
            \E\left[\theta^{\hat N_\tau}e^{-\lambda\tau (\theta-1)}
      1_{\{\hat T_1<\tau\}}
       \left((L^{-y}_ {\tau}(B)
       -L^{-y}_{\hat T_1}(B)
         \right) \right] \\
         &=&
         \theta
          \E\left[e^{\lambda\tau }
      1_{\{\hat N_\tau=1\}}
       \left((L^{-y}_ {\tau}(B)
       -L^{-y}_{\hat T_1}(B)
         \right) \right] +O(\theta^2).
         \end{eqnarray*}
We then have\begin{eqnarray*}
\bar{v}_{\lambda,\beta, \theta}(y)&=&\theta^{3/2}
              {v}_{\lambda,\beta}(y)+o(\theta^{3/2}),
              \end{eqnarray*}
              with
\begin{eqnarray*}
{v}_{\lambda,\beta}(y)&=&
 \sup_{\tau\in \T_{0,1}}
\E\left[ e^{\lambda\tau} 1_{\{\hat N_\tau=0\}}
\int_0^{\tau} f_{\lambda\beta}(
   y+ B_{s}
    )ds
     +\frac{\beta}{2} 
     e^{\lambda\tau} 1_{\{\hat N_\tau=1\}}
      \left(L^{-y}_ {\tau}(B)
       -L^{-y}_{\hat T_1}(B)
         \right) 
          \right].
\end{eqnarray*}
Finally, we obtain
 \begin{eqnarray*}
P(T-\theta,b(T)e^{a\sqrt{\theta}})-(K-b(T)e^{a\sqrt{\theta}})
&=&
\theta^\frac{3}{2}(\sigma b(T)\bar\delta e^\lambda) \vv_{\lambda,\beta}\left (\frac{a}{\sigma}\right )+o(\theta^\frac{3}{2}),
\end{eqnarray*}

\end{dema}

\subsection{Convergence rate of the critical price}\label{S2}
Thanks to the expansion given in Theorem~\ref{thjd1},  we  are now able to state the first main result of this paper.
\begin{theorem} \label{thjd2}
Under the hypothesis of the model and  $\bar d<0$, we have :
\\If $\nu\{\ln{K\over b(T)}\}=0$, then we have 
 $$\lim_{t\rightarrow T }\frac{b(T)-b(t)}{\sigma b(T)\sqrt{(T-t)}}= y_{0},$$
  with $ y_{0}=-\sup\{x\in \R\;;\;v_{0}(x)=\sup_{\tau \in \T_{0,1}}\E(\int_0^\tau (x+B_s)ds)=0\}$.\\
 If $\nu\{\ln{K\over b(T)}\}>0$, 
 we then  have 
 $$\lim_{t\rightarrow T }\frac{b(T)-b(t)}{\sigma b(T)\sqrt{(T-t)}}= {y_{\lambda,\beta}},$$
  with $ {y_{\lambda,\beta}}$ as defined in  Lemma \ref{v_alpha}, with 
$${{\lambda}}={\nu\{\ln{K\over b(T)}\} }\;,\; 
\beta={K\over b(T)\bar\delta  }
 \;\mbox{and}\quad \bar{\delta}=\delta+\int_{y>\ln(K/b(T))}\hspace{-11mm}e^y\nu(dy)   .$$ 
\end{theorem}
\begin{dema}\textbf{of Theorem \ref{thjd2}:\\}{
\medskip
According to  Theorem \ref{thjd1}, we have for all $a<0$,\[P(T-\theta,b(T)e^{a\sqrt{\theta}})=(K-b(T)e^{a\sqrt{\theta}})^++ C\theta^\frac{3}{2}\vv_{\lambda,\beta}(\frac{a}{\sigma})+o(\theta^\frac{3}{2}).\]
\noindent { \bf Lower bound for   $b(T)-b(t)$}
\\
Specifically, we have for all $a>-\sigma {y_{\lambda,\beta}}$, where ${y_{\lambda,\beta}}$ is defined by Lemma \ref{v_alpha},  $$\vv_{\lambda,\beta}(\frac{a}{\sigma})>0,$$ we thus obtain for $\theta$ close to  0, $$P(t,b(T)e^{a\sqrt\theta}) > (K-b(T)e^{a\sqrt\theta}),$$ and then $$\ln(b(T))+a\sqrt{\theta}>\ln(b(t)),$$
hence
 \begin{eqnarray*}
\frac{b(T)-b(t)}{b(t)\sqrt{\theta}}>-a.\\
\end{eqnarray*}
Noting that since $r>0$ we have $b(T)>0$, and by making $t$ tend to $T$ then  $a$ to $-\sigma {y_{\lambda,\beta}}$, we  obtain
$$\liminf_{t\rightarrow T}\frac{b(T)-b(t)}{b(T)\sqrt{T-t}}\geq\sigma {y_{\lambda,\beta}}.$$
\medskip

\noindent { \bf Upper bound for   $b(T)-b(t)$}
\\
Let's consider $a\leq - \sigma {y_{\lambda,\beta}}$, we have thus $\vv_{\lambda,\beta}(\frac{a}{\sigma})=0$ and consequently, $$P(t,b(T)e^{a\sqrt\theta}) - (K-b(T)e^{a\sqrt\theta})=g(\theta),$$ with $g(\theta)=o(\theta^{\frac{3}{2}})$.
\\In addition, we have for all $b(t)<x<K$,
$$P(t,x)-P(t,b(t))-(x-b(t))\frac{\partial P}{\partial x}(t,b(t))=\int_{b(t)}^x(u-b(t))\frac{\partial^2 P}{\partial x^2}(t,du), $$
since $\frac{\partial^2 P}{\partial x^2}(t,du)$ is a positive measure on $]0,+\infty[$. As the smooth-fit is satisfied, $\frac{\partial P}{\partial x}(t,b(t))=-1$ (see \cite{mikousmooth}), we have for all $b(t)<x<K$, 
\begin{eqnarray*}
P(t,x)-(K-x)&=&\int_{b(t)}^x(u-b(t))\frac{\partial^2 P}{\partial x^2}(t,du).\\
\end{eqnarray*}
Then, for $b(t)<x=b(T)e^{a\sqrt{\theta}}$, we have according to Lemma \ref{min_deriv2_P}, 
\begin{eqnarray*}
\frac{u^2\sigma^2}{2}\frac{\partial^2P}{\partial x^2}(t,u)&\geq& b(T)\bar\delta\left ((1-e^{a\sqrt{\theta}})-{{\lambda\beta}} \sqrt{\theta}\sigma\E\left ( B_1+\frac{a}{\sigma}\right )^+\right )+o(\sqrt{\theta})
\\
&\geq& b(T)\bar\delta\sqrt{\theta}\sigma\left (-\frac{a}{\sigma}-{{\lambda\beta}} \E\left ( B_1+\frac{a}{\sigma}\right )^+\right )+o(\sqrt{\theta})\\
\end{eqnarray*}
Hence 
$$
P(t,x)-(K-x)\geq [(x-b(t))^+]^2\left (\frac{C(-\frac{a}{\sigma})}{b(0)^2\sigma^2}\sqrt{\theta}+o(\sqrt{\theta})\right ),$$ where  $C(x)=x-{{\lambda\beta}}\E(B_1-x)^+$.
Due to lemma \ref{conjecture} and to the  continuity of $C(x)$, we have, for $\frac{-a}{\sigma}$ close enough to $ {y_{\lambda,\beta}}$, $C(-\frac{a}{\sigma})>0$. Moreover,
\begin{eqnarray*}
P(t,b(T)e^{a\sqrt\theta})-(K-b(T)e^{a\sqrt\theta})=g(\theta)=o(\theta^\frac{3}{2}).
\end{eqnarray*}
Therefore, for $\theta$  small enough, there exists a positive constant $A$ such that 
\begin{equation}\label{reutilisable}
[(b(T)e^{a\sqrt\theta}-b(t))^+]^2 \leq A b(0)^2\sigma^2 \frac{g(\theta)}{C(-\frac{a}{\sigma})\sqrt\theta}=o(\theta) 
\end{equation} 
$$(b(T)e^{a\sqrt\theta}-b(t))^+= o(\sqrt\theta),$$
and then, for $\theta$ small enough, 
$$\frac{b(T)-b(t)}{b(T)\sqrt\theta}\leq -a+ o(1).$$ 
Finally, by making $a$ tend to $-\sigma {y_{\lambda,\beta}}$, we obtain
$$\limsup_{t\rightarrow T} \frac{b(T)-b(t)}{b(T)\sqrt{T-t}}\leq \sigma {y_{\lambda,\beta}}.$$
}\end{dema}
\section{Limit case}

In this part, we consider the limit case where $\bar d=r-\delta-\int_{y>0}(e^y-1)\nu(dy)=0$, we then have  
\begin{theorem} \label{thjd6}
According to the model hypothesis, if $\bar d=0$, then, we  have  
$$\lim_{t\rightarrow T }\frac{K-b(t)}{\sigma K\sqrt{(T-t)|\ln(T-t)|}}=\sqrt{2}.$$
\end{theorem}

The method for proving Theorem \ref{thjd6} consists of analysing the behavior of  the European critical price $b_e(t)$ introduced in section 1, afterwards we prove that the behavior of the  critical price $b(t)$ is similar by controling the difference $b(t)-b_e(t)$.

Let us denote by \[{\alpha}(\theta)=\frac{\ln(\frac{K}{b_e(t)})-\mu\theta}{\sigma \sqrt\theta},\]
where $\mu=\gamma_0-\frac{\sigma^2}{2}=r-\delta-\int(e^y-1)\nu(dy)-\frac{\sigma^2}{2}$.
 \begin{prop}
\label{rate_b_e}
Under the model hypothesis, if $\bar d=0$, then we have\\
\begin{itemize}
\item[i)]${\alpha}(\theta)\sim \sqrt{2\ln(\frac{1}{\theta})}$
\item[ii)]$\lim_{\theta \rightarrow0}\frac{K-b_e(t)}{\sigma K\sqrt{|\theta\ln(\theta)|}}= \sqrt2$
\end{itemize}

\end{prop}
\pproof{of Proposition \ref{rate_b_e}}
{
 Since $b(t)\leq b_e(t)\leq K$ and $b(t)\rightarrow K$, we clearly have
 $$\sqrt{\theta}{\alpha}(\theta)\longrightarrow_{\theta\rightarrow 0}0.$$
We will first prove that  ${\alpha}(\theta)\longrightarrow_{\theta\rightarrow 0}+\infty$, or equivalently.
\begin{equation}\label{dl1}
\lim_{\theta\rightarrow 0} \frac{K-b_e(t)}{\sigma\sqrt{\theta}}=+\infty.
\end{equation}
We have
 \begin{equation}\label{prix euro}K-b_e(t)=e^{-r\theta}\E\left [\left (K-b_e(t)e^{\tilde{X}_\theta}\right )^+\right ]\end{equation}
Therefore 
 \begin{eqnarray*}
 \frac{K-b_e(t)}{\sqrt{\theta}}& =&e^{-r\theta}\E\left [\left (\frac{K-b_e(t)}{\sqrt{\theta}}+b_e(t) \frac{1-e^{\tilde{X}_\theta}}{\sqrt{\theta}}\right )^+\right ]\\
 & =&e^{-r\theta}\E\left [\left (\frac{K-b_e(t)}{\sqrt{\theta}}+b_e(t) \frac{1-e^{\sigma\sqrt\theta{B}_1+\mu\theta +Z_\theta}}{\sqrt{\theta}}\right )^+\right ]\\
 \end{eqnarray*}
 Now, if we notice that $\frac{1-e^{\sigma\sqrt\theta{B}_1+\mu\theta +Z_\theta}}{\sqrt{\theta}}\longrightarrow_{\theta \rightarrow0}^{p.s}-\sigma B_1$, we have by Fatou lemma
 \begin{eqnarray*}
 \liminf_{\theta\rightarrow 0}\frac{K-b_e(t)}{\sqrt{\theta}} &\geq& \E\left [\left (\liminf_{\theta\rightarrow 0}\frac{K-b_e(t)}{\sqrt{\theta}}-\sigma K B_1\right )^+\right ]\\
&=& \liminf_{\theta\rightarrow 0}\frac{K-b_e(t)}{\sqrt{\theta}}+\E\left [\left (\sigma K B_1-\liminf_{\theta\rightarrow 0}\frac{K-b_e(t)}{\sqrt{\theta}}\right )^+\right ]
  \end{eqnarray*}
which is equivalent to
\[\E\left [\left (\sigma K B_1-\liminf_{\theta\rightarrow 0}\frac{K-b_e(t)}{\sqrt{\theta}}\right )^+\right ]\leq 0.\] 
This gives  (\ref{dl1}) which yields the  wanted result.

i) We now rewrite  equation (\ref{prix euro}) to obtain
\[K-b_e(t)=e^{-r\theta}K-b_e(t)e^{-\delta \theta}+e^{-r\theta}\E\left [\left (b_e(t)e^{\tilde{X}_\theta}-K\right )^+\right ],\]
 therefore
\begin{equation}\label{Ba}
 e^{-r\theta}\E\left [\left (e^{\tilde{X}_\theta}-e^{\ln(\frac{K}{b_e(t)})}\right )^+\right ]=\frac{K}{b_e(t)}(1-e^{-r\theta})-(1-e^{-\delta\theta}) .
 \end{equation}
We will give an expansion  for each side of the equation. For the left hand side of the equation, we have 
\begin{eqnarray*}
\lefteqn{e^{-r\theta}\E\left [\left (e^{\tilde{X}_\theta}-e^{\ln(\frac{K}{b_e(t)})}\right )^+\right ]}\\
&= &e^{-r\theta +\theta \mu+\sigma {\alpha}(\theta)\sqrt{\theta}}\E\left [\left (e^{\sigma \sqrt \theta B_1+Z_\theta -\sigma{\alpha}(\theta) \sqrt{\theta}}-1\right )^+\right ]\\
&=&e^{-r\theta +\theta \mu+\sigma {\alpha}(\theta)\sqrt{\theta}}\E\left [\left (U_\theta e^{Z_\theta }-1\right )^+\right ],\\
\end{eqnarray*}
where $U_\theta=e^{\sigma \sqrt \theta B_1 -\sigma{\alpha}(\theta) \sqrt{\theta}}$.
Since the process $Z_t$ is independent of  $U_\theta$, we can write 
\begin{eqnarray*}
\lefteqn{\E\left [\left (U_\theta e^{Z_\theta }-1\right )^+|U_\theta\right ]}\\
&=&\left (U_\theta -1\right )^++\E\left [\int_0^\theta ds\int \left (\left (U_\theta e^{Z_s+y}-1\right )^+-\left (U_\theta e^{Z_s }-1\right )^+\right )\nu(dy)|U_\theta\right  ]\\
&=&\left (U_\theta -1\right )^++\int_0^\theta ds\int \left (\left (U_\theta e^{y}-1\right )^+-\left (U_\theta -1\right )^+\right )\nu(dy)+U_\theta O(\theta^2),
\end{eqnarray*}
where $O(\theta^2)$ is detrministic. Indeed,
\begin{eqnarray*} 
\lefteqn{\left |\E\left [\int_0^\theta ds\int \left (\left (U_\theta e^ye^{Z_s}-1\right )^+-\left (U_\theta e^y-1\right )^+\right )\nu(dy)|U_\theta\right  ]\right |}\\
&\leq &U_\theta \int e^y\nu(du)\int_0^\theta  \E\left | e^{Z_s}-1\right|ds=U_\theta O(\theta^2).\\
\end{eqnarray*}
 Taking the expectation, we thus obtain,
\begin{eqnarray*}
\lefteqn{\E\left [\left (U_\theta e^{Z_\theta }-1\right )^+\right ]}\\
&&=\E\left [\left (U_\theta -1\right )^+\right ]+\theta\int\E\left [\left (U_\theta e^{y}-1\right )^+\right ]\nu(dy)-\nu(\R)\theta\E\left [\left (U_\theta -1\right )^+\right ]+O(\theta^2).
\end{eqnarray*}
Since $\alpha(\theta)\rightarrow \infty$, we have like in \cite{damien} $$\E\left [\left (U_\theta -1\right )^+\right ]\sim \sigma\sqrt{\theta}\E\left (B_1 -{\alpha}(\theta)\right )^+=o(\sqrt\theta),$$
 then
\begin{eqnarray*}
\E\left [\left (U_\theta e^{Z_\theta }-1\right )^+\right ]
&=&\E\left [\left (U_\theta -1\right )^+\right ]+\theta\int\E\left [\left (U_\theta e^{y}-1\right )^+\right ]\nu(dy)+o(\theta^\frac{3}{2}).
\end{eqnarray*}
We recall that $U_\theta=e^{\sigma \sqrt \theta B_1 -\sigma{\alpha}(\theta) \sqrt{\theta}}$, then  
\begin{eqnarray*}
\lefteqn{\E\left [\left (U_\theta e^{y}-1\right )^+\right ]-\left (e^{y -\sigma{\alpha}(\theta) \sqrt{\theta}} -1\right )^+}\\
&\leq&e^{y -\sigma{\alpha}(\theta) \sqrt{\theta}}\E\left |e^{\sigma \sqrt \theta B_1}-1\right | = e^yO(\sqrt\theta).
\end{eqnarray*}
Hence, 
\begin{eqnarray*}
\lefteqn{\E\left [\left (U_\theta e^{Z_\theta }-1\right )^+\right ]}\\
&=&\E\left [\left (U_\theta -1\right )^+\right ]+\theta\int\left ( e^{y-\sigma{\alpha}(\theta) \sqrt{\theta}}-1\right )^+\nu(dy)+O(\theta^\frac{3}{2})\\
&=&\E\left [\left (U_\theta -1\right )^+\right ]+\theta\int_{y>0}\left ( e^{y-\sigma{\alpha}(\theta) \sqrt{\theta}}-1\right )\nu(dy)\\
&&-\theta\int_{0<y<\sigma{\alpha}(\theta) \sqrt{\theta}}\left ( e^{y-\sigma{\alpha}(\theta) \sqrt{\theta}}-1\right )\nu(dy)+O(\theta^\frac{3}{2}).
\end{eqnarray*}
 Since  $(1-e^{-x})\leq x$, we  then have
\begin{eqnarray*}
\left |\int_{\left (0,\sigma{\alpha}(\theta) \sqrt{\theta}\right )} ( e^{y-\sigma{\alpha}(\theta) \sqrt{\theta}}-1 )\nu(dy)\right |
&\leq& \underbrace{\overbrace{\nu\{0<y<\sigma{\alpha}(\theta) \sqrt{\theta}\}}^{\longrightarrow_{\theta \rightarrow
0}0} \sigma{\alpha}(\theta) \sqrt{\theta}}_{=o({\alpha}(\theta)\theta^\frac{1}{2}) },
\end{eqnarray*}
and noticing that $\theta^\frac{3}{2}=o(\alpha(\theta)\theta^\frac{3}{2}) $, we obtain
\begin{eqnarray*}
\E\left (U_\theta e^{Z_\theta }\hspace{-2mm}-1\right )^+
\hspace{-2mm}=\E\left (U_\theta \hspace{-1mm}-1\right )^+\hspace{-2mm} +\theta \hspace{-2mm} \int\hspace{-2mm} \left ( e^{y}-1\right )^+\hspace{-2mm} \nu(dy)
 -{\alpha}(\theta) \theta^\frac{3}{2}\sigma \hspace{-2mm} \int_{y>0}\hspace{-4mm} e^y\nu(dy)+o({\alpha}(\theta) \theta^\frac{3}{2}).\\
\end{eqnarray*}
The left hand side of  equation (\ref{Ba}) becomes
\begin{eqnarray}
\lefteqn{e^{-r\theta}\E\left [\left (e^{\tilde{X}_\theta}-e^{\ln(\frac{K}{b_e(t)})}\right )^+\right ]
=e^{-r\theta +\theta \mu+\sigma {\alpha}(\theta)\sqrt{\theta}}\E\left [\left (U_\theta e^{Z_\theta }-1\right )^+\right ]}\nonumber\\
&=&e^{-r\theta +\theta \mu+\sigma {\alpha}(\theta)\sqrt{\theta}}\E\left [\left (U_\theta -1\right )^+\right ]\nonumber\\
&+&\left (1+\sigma {\alpha}(\theta)\sqrt{\theta}+o({\alpha}(\theta) \sqrt{\theta})\right )\left (\theta\int\left ( e^{y}-1\right )^+\nu(dy)\right. 
\left .-{\alpha}(\theta) \theta^\frac{3}{2}\sigma\int_{y>0}\hspace{-3MM}e^y\nu(dy)+ o({\alpha}(\theta) \theta^\frac{3}{2})\right )\nonumber\\
&=&e^{-r\theta +\theta \mu+\sigma {\alpha}(\theta)\sqrt{\theta}}\E\left [\left (U_\theta -1\right )^+\right ]+\theta\int(e^y-1)^+\nu(dy)
-\nu(\R^+)\sigma{\alpha}(\theta) \theta^\frac{3}{2}+o({\alpha}(\theta) \theta^\frac{3}{2}).\nonumber\\\label{eqiv1}
\end{eqnarray}
Besides, the right hand side of  (\ref{Ba})
\begin{eqnarray}
\lefteqn{\frac{K}{b_e(t)}(1-e^{-r\theta})-(1-e^{-\delta\theta})=
e^{\sigma\sqrt\theta{\alpha}(\theta)+\mu\theta} r\theta-\delta \theta +O(\theta^2)}\nonumber\\
&= &(r-\delta)\theta + r\sigma\theta^\frac{3}{2}{\alpha}(\theta)+o(\theta^\frac{3}{2}{\alpha}(\theta))\nonumber\\
&= & \left (\hspace{-0mm}\int\hspace{0mm}(e^y-1)^+\nu(dy) \hspace{-0mm}\right )\theta + r\sigma {\alpha}(\theta)\theta^\frac{3}{2} +o(\theta^\frac{3}{2}{\alpha}(\theta)).\label{eqiv2}
\end{eqnarray}
Thanks to (\ref{eqiv1}) and (\ref{eqiv2}),  equation (\ref{Ba}) becomes, 
$$e^{-r\theta +\theta \mu+\sigma {\alpha}(\theta)\sqrt{\theta}}\E\left [\left (U_\theta -1\right )^+\right ]=\sigma\left (r+\nu(\R^+)\right ) {\alpha}(\theta)\theta^\frac{3}{2} +o(\theta^\frac{3}{2}{\alpha}(\theta)).
$$
Hence, 
$$\E\left [\left (U_\theta -1\right )^+\right ]\sim \sigma\left (r+\nu(\R^+)\right ) {\alpha}(\theta)\theta^\frac{3}{2} .
$$
As explained above, thanks to proposition 2.1 in \cite{villeneuve}, we have
$$\E\left [\left (U_\theta -1\right )^+\right ]\sim \sigma\sqrt\theta\E(B_1-{\alpha}(\theta))^+ \sim \frac{\sigma \sqrt\theta}{\sqrt{2\pi}{\alpha}^2(\theta)e^{\frac{{\alpha}^2(\theta)}{2}}} .
$$
Thus, we have 
\begin{equation}\label{E(b-a)+}
\frac{1}{\sqrt{2\pi}{\alpha}^2 (\theta)e^{\frac{{\alpha}^2(\theta)}{2}}}\sim \left (r+\nu(\R^+)\right )\theta{\alpha}(\theta),
\end{equation}
 hence \begin{equation}\label{eqiv alpha}\alpha(\theta)\sim \sqrt{2\ln(\frac{1}{\theta})}.\end{equation}

ii) Since $\frac{K-b_e(t)}{K\sigma\sqrt\theta}\sim{\alpha}(\theta)$, we obtain $$\frac{K-b_e(t)}{\sigma K}\sim\sqrt{2\theta\ln(\frac{1}{\theta})}.$$
}
To compare the behaviors of $b(t)$ and $b_e(t)$, we have to control the difference between them.
\begin{prop}\label{b-b_e}
According to the model hypothesis, if $\bar d=0$, then there exists $C>0$ such that
\[0\leq\frac{b_e(t)- b(t)}{\sqrt{T-t}}\leq C.\]
\end{prop}
Before proving Proposition \ref{b-b_e}, we need to prove the non decreasing of $b_e(t)$ near maturity which is the purpose of this following lemma
\begin{lemma}\label{be_increase}
The critical European put price, $b_e(t)$, is differentiable on $\left (0,T\right )$ and for  $t$ close to $T$, we have\[b'_e(t)\geq0.\]
\end{lemma}
\pproof{ of Lemma \ref{be_increase}}{
We recall that $F$ is the function defined by $F(t,x)= P_e(t,x)-(K-x)$,  $F$ is $\mathcal{C}^1$  on $\left (0, T\right )\times \left (0, K\right )$ and  satisfies $\frac{\partial F}{\partial x}(t,x)=\frac{\partial P_e}{\partial x}(t,x)+1>0  $.
Due to its definition, $b_e(t)$ satisfies the following equation,
$P_e(t,b_e(t))-(K-b_e(t))=0$. Then, thanks to the implicit function theorem, $b_e(t)$ is differentiable on $\left (0, T\right )$ and \[b_e'(t)=-\frac{\frac{\partial F}{\partial t}(t,b_e(t))}{\frac{\partial F}{\partial x}(t,b_e(t))}=-\frac{\frac{\partial P_e}{\partial t}(t,b_e(t))}{\frac{\partial P_e}{\partial x}(t,b_e(t))+1},\]
which means  that 
\[-b_e'(t)\frac{\partial P_e}{\partial t}(t,b_e(t))\leq0.\]
We will study the sign of $\frac{\partial P_e}{\partial t}(t,b_e(t))$ instead of that of $b_e'(t)$.\\
The European put price satisfies the following equation
\begin{eqnarray*} 
\lefteqn{\frac{\partial P_e}{\partial t}(t,b_e(t)) }\\
&=&rP_e(t,b_e(t))-\frac{\sigma^2 b_e(t)^2}{2}\frac{\partial^2P_e}{\partial x^2}(t,b_e(t))-(r-\delta)b_e(t)\frac{\partial P_e}{\partial x}(t,b_e(t))\\
&&-\int \left [P_e(t,b_e(t)e^y)-P_e(t,b_e(t))-b_e(t)(e^y-1)\frac{\partial P_e}{\partial x}(t,b_e(t))\right ]\nu(dy)\\
&=&r(K-b_e(t))-\frac{\sigma^2 b_e(t)^2}{2}\frac{\partial^2P_e}{\partial x^2}(t,b_e(t))-\underbrace{ \left (\hspace{-1mm}r\hspace{-1mm}-\hspace{-1mm} \delta\hspace{-1mm}-\hspace{-2mm}\int\hspace{-1mm}(e^y-1)^+\hspace{-1mm}\nu(dy)\hspace{-1mm}\right )}_{\bar d=0}b_e(t)\frac{\partial P_e}{\partial x}(t,b_e(t))\\
&&-\int_{y>0} P_e(t,b_e(t)e^y)\nu(dy)+\nu(\R_+)P_e(t,b_e(t))\\
&&-\int_{y<0} \hspace{-1,5mm}\left [P_e(t,b_e(t)e^y)-P_e(t,b_e(t))-b_e(t)(e^y-1)\frac{\partial P_e}{\partial x}(t,b_e(t))\right ]\nu(dy)\\
\end{eqnarray*}
Since $P_e(t,.)$ is a non negative  convex function, we have $\int_{y>0} \left [P_e(t,b_e(t)e^y)\right ]\nu(dy)\geq 0$ and
\[\int_{y<0} \hspace{-1,5mm}\left [P_e(t,b_e(t)e^y)-P_e(t,b_e(t))-b_e(t)(e^y-1)\frac{\partial P_e}{\partial x}(t,b_e(t))\right ]\nu(dy)\geq0,\]
so that
\begin{eqnarray*} 
\frac{\partial P_e}{\partial t}(t,b_e(t))
&\leq &\left (r+\nu(\R^+)\right )(K-b_e(t))-\frac{\sigma^2 b_e(t)^2}{2}\frac{\partial^2P_e}{\partial x^2}(t,b_e(t)).\\
\end{eqnarray*}
Thanks to lemma \ref{be_increase}, we have an equivalent for $(K-b_e(t))$. Now, let's have a look  at the estimate of $\frac{\partial^2P_e}{\partial x^2}(t,b_e(t))$  near $T$. We have 
\begin{eqnarray*} 
\frac{\partial P_e}{\partial x} (t,x)&=&-e^{-r(T-t)}\E\left [e^{\tilde{X}_{T-t}}1_{\{K-xe^{\tilde{X}_{T-t}}>0\}}\right ]
\\&=&
-e^{-r(T-t)}\int_{-\infty}^{\ln(\frac{K}{x})}e^u p_{\tilde{X}_{T-t}}(u) du,
\end{eqnarray*}
where $p_{X}$ denotes the density of $X$ and $\tilde{X}_{t}=\mu(t)+\sigma B_{t}+Z_{t}$. 
Then, we have 
\begin{eqnarray*} 
\frac{\partial^2 P_e}{\partial x^2} (t,x)&=&
e^{-r(T-t)}\frac{K}{x^2}\ p_{\tilde{X}_{T-t}}\left (\ln(\frac{K}{x})\right )\\
&\geq& e^{-r(T-t)}\frac{K}{x^2} p_{\mu(T-t)+\sigma B_{T-t}}\left (\ln(\frac{K}{x})\right )\PP(T_1>\theta)\\
									  &= &e^{-r\theta}\frac{K}{x^2\sigma \sqrt{2\pi \theta}}e^{\frac{-1}{2}\left (\frac{\ln\left (K/x\right )-\mu \theta}{\sigma\sqrt{\theta}}\right )^2}\PP(T_1>\theta).
\end{eqnarray*}
Then,
\begin{eqnarray*} 
\frac{\partial P_e}{\partial t}(t,b_e(t)) 
&\leq &\left (r+\nu(\R^+)\right )(K-b_e(t))-e^{-r\theta}\frac{\sigma K}{2\sqrt{2\pi \theta}}e^{-\frac{\alpha(\theta)^2}{2}}\PP(T_1>\theta).
\end{eqnarray*}
We can easily check that $K-b_e(t)\sim\sigma K\sqrt{\theta}\alpha(\theta)=o({\alpha}^3(\theta) \sqrt\theta)$, and we recall  the equivalency (\ref{E(b-a)+})
\[\frac{1}{\sqrt{2\pi}{\alpha}^2 (\theta)e^{\frac{{\alpha}^2(\theta)}{2}}}\sim \left (r+\nu(\R^+)\right )\theta{\alpha}(\theta),\]
which yields $$e^{-r\theta}\frac{\sigma K}{2\sqrt{2\pi \theta}}e^{-\frac{{\alpha}(\theta)^2}{2}}\PP(T_1>\theta)\sim\sigma K\frac{ e^{-\frac{{\alpha}^2(\theta)}{2}}}{2\sqrt{2\pi \theta}}\sim\frac{\sigma K}{2} \left (r+\nu(R^+)\right ){\alpha}^3(\theta)\sqrt\theta.$$
Then, we have, for $\theta$ small enough
\begin{eqnarray*}
\frac{\partial P_e}{\partial t}(t,b_e(t))
&\leq&-\frac{\sigma K}{2} \left (r+\nu(R^+)\right ){\alpha}^3(\theta)\sqrt\theta +o({\alpha}^3(\theta)\sqrt\theta)<0,
\end{eqnarray*} 
which proves that $b_e'(t)$ is a non decreasing function for $t$ close to $T$.
}
We are now in a position to prove  Proposition \ref{b-b_e}.

\pproof{ of Proposition \ref{b-b_e}}{
An expansion of   $P(t,x)$ around $(t,b(t))$ gives
$$P(t,x)-P(t,b(t))-(x-b(t))\frac{\partial P}{\partial x}(t,b(t))=\int_{b(t)}^x(u-b(t))\frac{\partial^2 P}{\partial x^2}(t,du), $$
and thanks to the smooth-fit which is satisfied at $b(t)$, we obtain
\begin{eqnarray*}
P(t,x)-(K-x)&\geq& \frac{\left (x-b(t)\right )^2}{2}\inf_{b(t)\leq u \leq x} \frac{\partial^2P}{\partial x^2}(t,u).
\end{eqnarray*}
First, we are going to give, as in Lemma  \ref{min_deriv2_P}, a lower bound for $\underset{b(t)\leq u \leq b_e(t)} {\inf}\frac{u^2\sigma^2}{2}\frac{\partial^2P}{\partial x^2}(t,u)$.
The variational inequality gives, for    $u\in \left(b(t),K\right)$,
\begin{eqnarray*}
\lefteqn{\frac{u^2\sigma^2}{2}\frac{\partial^2P}{\partial x^2}(t,u)}
\\&\geq& rP(t,u)-(r-\delta)u\frac{\partial P}{\partial x}(t,u)-\hspace{-2mm}\int \hspace{-2mm}\left (\hspace{-1mm}P(t,ue^y)-P(t,u)-u(e^y-1)\frac{\partial P}{\partial x}(t,u)\hspace{-1mm}\right)\hspace{-1mm}\nu(dy)
\\&\geq& r(K-u)-\hspace{-1mm}\left (\hspace{-1mm}r-\delta-\int_{y>0}\hspace{-2mm}\hspace{-2mm}(e^y-1)\hspace{-0mm}\nu(dy)\hspace{-1mm}\right)u\frac{\partial P}{\partial x}(t,u)-\hspace{-1mm}\int_{y>0} \hspace{-4mm}P(t,ue^y)-P(t,u)\hspace{-0mm}\nu(dy)
\\&&-\int_{y<0} \hspace{-2mm}\left (\hspace{-1mm}P(t,ue^y)-(K-u)-u(e^y-1)\frac{\partial P}{\partial x}(t,u)\hspace{-1mm}\right)\hspace{-1mm}\nu(dy)
.\\
\end{eqnarray*}
Since $P(t,.)$ is non increasing and $\bar d=0$, we obtain 
\begin{eqnarray*}
\lefteqn{
\frac{u^2\sigma^2}{2}\frac{\partial^2P}{\partial x^2}(t,u) 
\geq r(K-u)-\int_{y<0} \hspace{-2mm}\left (\hspace{-1mm}P(t,ue^y)-(K-u)-u(e^y-1)\frac{\partial P}{\partial x}(t,u)\hspace{-1mm}\right)\hspace{-1mm}\nu(dy)
}\\
&=&r(K\hspace{-1mm}-\hspace{-1mm}u)\hspace{-1mm}-\hspace{-2mm}\int_{y<0} \hspace{-5mm} P(t,ue^y)\hspace{-1mm}-\hspace{-1mm}(K-ue^y)\nu(dy)-u\left(\frac{\partial P}{\partial x}(t,u)+1\right)\left(\int_{y<0}\hspace{-4mm}(1-e^y)\nu(dy)\right)
.
\end{eqnarray*}
Thanks to the convexity of $P$, $\frac{\partial P}{\partial x}(t,u)$ is non decreasing and $\frac{\partial P}{\partial x}(t,u)\geq-1$ . We then have, for all $t<T$,
\begin{eqnarray*}
\inf_{b(t)\leq u\leq b_e(t)}
\frac{u^2\sigma^2}{2}\frac{\partial^2P}{\partial x^2}(t,u) 
&\geq& r(K\hspace{-1mm}-\hspace{-1mm}b_e(t))\hspace{-1mm}-\hspace{-2mm}\int_{y<0} \hspace{-5mm} P(t,b_e(t)e^y)\hspace{-1mm}-\hspace{-1mm}(K-b_e(t)e^y)\nu(dy)
\\&&-b_e(t)\left(\frac{\partial P}{\partial x}(t,b_e(t))+1\right)\left(\int_{y<0}\hspace{-4mm}(1-e^y)\nu(dy)\right)\\
&\geq& r(K\hspace{-1mm}-\hspace{-1mm}b_e(t))\hspace{-1mm}-\hspace{-2mm}\int_{y<0} \hspace{-5mm} P_e(t,b_e(t)e^y)\hspace{-1mm}-\hspace{-1mm}(K-b_e(t)e^y)\nu(dy)+o(\sqrt\theta)
.
\end{eqnarray*}
We obtained the last inquality, using the estimate of $e(\theta,x)=O(\theta)$ and $\frac{\partial P}{\partial x}(t,x)+1=o(\sqrt\theta)$ (see Lemma \ref{deriv_prime}).  Since $y<0$, we also have $P_e(t,b_e(t)e^y)\hspace{-1mm}-\hspace{-1mm}(K-b_e(t)e^y)\leq0$, thus
\begin{eqnarray*}
\inf_{b(t)\leq u\leq b_e(t)}
\frac{u^2\sigma^2}{2}\frac{\partial^2P}{\partial x^2}(t,u) 
&\geq& r(K\hspace{-1mm}-\hspace{-1mm}b_e(t))+o(\sqrt\theta)
.
\end{eqnarray*}
 Besides,  for $\theta$ small enough, we have  $\sqrt\theta \leq K-b_e(t)$, then we obtain
\begin{eqnarray*}
P(t,b_e(t))-(K-b_e(t))
&\geq&  \frac{[(b_e(t)-b(t))^+]^2}{b_e^2(t)\sigma^2}r(K-b_e(t))(1+o(1)).
\end{eqnarray*}
Furthermore,
\begin{eqnarray*}
\lefteqn{P(t,b_e(t))-(K-b_e(t))=e(\theta ,b_e(t))}\\
&=&\E \left\{\hspace{-1MM}
\int_0^\theta \hspace{-3MM}e^{-rs}\hspace{-1MM}
\left(
 \hspace{-1MM}rK\hspace{-1MM}-\hspace{-1MM}\delta S^{b_e(t)}_s\hspace{-1MM}-\hspace{-2MM}\int_{y>0} \hspace{-6MM} P(t+s,S^{b_e(t)}_se^y)\hspace{-1MM}- \hspace{-1MM}\left(\hspace{-1MM}K\hspace{-1MM} -\hspace{-1MM}S^{b_e(t)}_s e^y\hspace{0MM}\right)\hspace{-1MM}\nu(dy) \hspace{-1MM}
\right) \hspace{-1MM}
 1_{\{S^{b_e(t)}_s<b(t+s)\}} \hspace{-1MM}ds\hspace{-1MM}
 \right\}\\
 &\leq&\E \left\{
\int_0^\theta e^{-rs}
\left(rK-\delta S^{b_e(t)}_s-\int_{y>0} \left(S^{b_e(t)}_s e^y-K\right)^+\nu(dy) \right)  1_{\{S^{b_e(t)}_s<b(t+s)\}} ds
 \right\}
\end{eqnarray*}
Since $\delta=r-\int_{y>0}(e^y-1)\nu(dy)$, we have 
\begin{eqnarray*}
\lefteqn{0 \leq\left (rK-\delta x-\int_{y>0} \left(x e^y-K\right)^+\nu(dy)\right )1_{\{x<b(t+s)\}}}\\
&\leq&\left( r(K- x)-\int_{y>0} \left(x e^y-K\right)^+-\left (x e^y-x\right )^+\nu(dy)\right) 1_{\{x<K\}}\\
&\leq &\left (r+\nu(\R^+)\right )\left(K-x\right)^+,
\end{eqnarray*}
 thus,
\begin{eqnarray*}
e(\theta ,b_e(t))
&\leq& \left (r+\nu(\R^+)\right )\E \left\{
\int_0^\theta e^{-rs}
\left(
 K-S_s^{b_e(t)}\right)^+ ds
 \right\}\\
&=& \left (r+\nu(\R^+)\right )
\int_0^\theta
 P_e(T-s,b_e(t)) ds\\
&=& \left (r+\nu(\R^+)\right )
\int_0^\theta
 P_e(t+u,b_e(t)) du.
\end{eqnarray*}
And as we saw in lemma \ref{be_increase}, near $T$, $b_e(t)$ is non-decreasing, then $b_e(t)\leq b_e(t+u)$. Due to the non-decreasing  of $P_e(t,x)- (K-x)$ on $x$, we thus have $$P_e(t+u,b_e(t))\leq K-b_e(t),$$ 
In conclusion, we have 
$$
e(\theta ,b_e(t))\leq\left (r+\nu(\R^+)\right )\theta(K-b_e(t))
$$
and
$$
e(\theta ,b_e(t))\geq\frac{[b_e(t)-b(t)]^2}{b^2(t)\sigma^2}\bar\delta (K-b_e(t)).
$$
Which gives the wanted result:
There  exists a constant $C$ such that
$$
\frac{b_e(t)-b(t)}{\sqrt{\theta}}\leq C.
$$
}
\section{Appendix 1: Proofs of lemmas}

\begin{dema}\textbf { of Lemma \ref{deriv_prime}:\\}{
According to  the early exercise premium  formula, we have, 
\begin{eqnarray*}
P(t,x)=P_e(t,x)+e(T-t,x)
\end{eqnarray*} 
and
\begin{eqnarray*}
&&e(\theta ,x)=\E \left\{
\int_0^\theta e^{-rs}
\Phi(t+s,xS^1_s)
 1_{\{xS_s^1<b(t+s)\}}  ds
 \right\},
\end{eqnarray*}
with $$\Phi(t,x)=rK-\delta x-\int_{y>0} P(t,xe^y)-\left(K-x e^y\right)\nu(dy).$$
Notice that $\Phi $ is a continuous function and   $\|\Phi'_x\|_\infty\leq \delta+\int_{y>0}e^y\nu(dy)$.
 \\1)  It is obvious that    $0 \leq e(\theta,x)\leq \theta rK=O(\theta)$, since $0\leq\Phi(t,x) 1_{\{x<b(t+s)\}}\leq rK$.
 \\2)
For all random variable  $X$, we denote by  
 $p_{X}$ its density,
we thus  have  for all fixed $s\in[0,\theta]$,  
\begin{eqnarray}
p_{-\tilde{X}_s}(x)&=&p_{-\mu s-\sigma B_s}*p_{-Z_s}(x)=\frac{1}{\sqrt{s}}\frac{1}{\sigma\sqrt{2\pi}}\int e^{-\frac{(-x+\mu s-u)^2}{2\sigma^2s}}p_{-Z_s}(u) du\nonumber\\
&\leq& C^{te}\frac{1}{\sqrt{s}}.\label{maj densit}
\end{eqnarray}
We can state
 \begin{eqnarray}
 \frac{\partial e}{\partial x} (\theta,x)
 &=&   \E \left\{
\int_0^\theta e^{-rs}S^1_s
\Phi'_x(t+s,xS^1_s)
 1_{\{xS_s^1<b(t+s)\}}  ds\right\}
 \nonumber\\
 &&   -\int_0^\theta \frac{\Phi(t+s,b(t+s))}{x}p_{-\tilde{X}_s}\left( \ln(\frac{x}{b(t+s)})\right )ds. \label{exist deriv prime}
 \end{eqnarray}
Then, we have 
\begin{eqnarray*}
\left |\frac{\partial e}{\partial x} (\theta,x)\right |
&\leq &  \left | \E \left\{
\int_0^\theta e^{-rs}S^1_s
\Phi'_x(t+s,xS^1_s)
 1_{\{xS_s^1<b(t+s)\}}  ds\right\}\right |\\
 &&   +\left |\int_0^\theta \frac{\Phi(t+s,b(t+s))}{x}p_{-\tilde{X}_s}\left( \ln(\frac{x}{b(t+s)})\right )ds\right |\\
&\leq & \|\Phi'_x\|_\infty \frac{b(T)}{x}\theta +  \left |\int_0^\theta \frac{\Phi(t+s,b(t+s))}{x}p_{-\tilde{X}_s}\left( \ln(\frac{x}{b(t+s)})\right )ds\right |.
\end{eqnarray*} 
 According to  inequality (\ref{maj densit})
, we also have
\begin{eqnarray*}
\lefteqn{\left |\int_0^\theta \frac{\Phi(t+s,b(t+s))}{x}p_{-\tilde{X}_s}\left( \ln(\frac{x}{b(t+s)})\right )ds\right |}\\
&&\leq C^{te}\left |\int_0^\theta \frac{\Phi(t+s,b(t+s))}{x\sqrt{s}}ds\right |\\
&&=\frac{C^{te}}{x} \theta \left |\int_0^1 \frac{\Phi(t+\theta u,b(t+\theta u))}{\sqrt{\theta u}}du\right |\\
&&\leq \frac{C^{te}}{x}\sqrt{\theta} \sup_{t\leq u\leq t+ \theta}\left |\Phi(u,b( u))\right |\int_0^1 \frac{1}{\sqrt{u}}du\\
&&\leq \frac{C^{te}}{x}\sqrt{\theta} \sup_{T-\theta\leq u\leq T}\left |\Phi( u,b( u))\right |\int_0^1 \frac{1}{\sqrt{u}}du.
\end{eqnarray*}
However, thanks to the continuity of $b(u)$ and of $\Phi( t,x)$, we have 
$\underset{\theta\rightarrow 0}{\lim}\underset{T-\theta\leq u\leq T}{\sup}\left |\Phi( u,b( u))\right |=\left |\Phi( T,b( T))\right |=0$.
Therefore, we conclude  that 
$\left |\frac{\partial e}{\partial x} (\theta,x)\right |=\frac{1}{x}o(\sqrt{\theta})$.
\\3)
Using the previous point, we have $\left|\frac{\partial e}{\partial x} (\theta,x)\right|=o(\sqrt\theta)$, then for all $x\leq b_e(t)\wedge b(T)$ and $\theta$ small enough , we have 
\begin{eqnarray*}
\lefteqn{
0\leq 1+\frac{\partial P}{\partial x}(t,x)
\leq  \left(1+\frac{\partial P_e}{\partial x}(t,x)\right)+o(\sqrt\theta)
}\\&=&1-\E\left[(e^{X_\theta}) 1_{\{{X_\theta}<\ln\frac{K}{x}\}}\right]+o(\sqrt\theta)
\\&\leq&1-\E\left[(e^{X_\theta}) 1_{\{{X_\theta}<\ln\frac{K}{b_e(t)\wedge b(T)}\}}\right]+o(\sqrt\theta).
\end{eqnarray*}\\
If $b(T)=K$, then for $x\leq b_e(t)$, 
\begin{eqnarray*}
1-\E\left[(e^{X_\theta}) 1_{\{{X_\theta}<\ln\frac{K}{b_e(t)\wedge b(T)}\}}\right]
&=&\PP({B_\theta}\geq\sqrt\theta\alpha(\theta))-\E\left({\sigma B_\theta} 1_{\{{B_\theta}<\sqrt\theta\alpha(\theta)\}}\right)+o(\sqrt\theta)
.
\end{eqnarray*}
Since $\alpha(\theta)\underset{\theta \rightarrow0}{\longrightarrow}\infty$, we have $$\left|\E\left({\sigma B_\theta} 1_{\{{B_\theta}<\sqrt\theta\alpha(\theta)\}}\right)\right|=\sigma\sqrt\theta\E\left({ B_1} 1_{\{B_1<\alpha(\theta)\}}\right)=o(\sqrt\theta), $$
and using  equivalencies (\ref{E(b-a)+}) and (\ref{eqiv alpha}), we also have
$$\PP\left({B_1}\geq\alpha(\theta)\right)\leq \frac{e^{-\frac{\alpha^2(\theta)}{2}}}{\alpha(\theta)}\leq C \theta \alpha^2(\theta)=O(\theta|\ln\theta|)=o(\sqrt\theta).$$\\
\\
If $b(T)<K$, then for $\theta$ small enough $b(T)<b_e(t)$ and 
\begin{eqnarray*}
\lefteqn{
1-\E\left[(e^{X_\theta}) 1_{\{{X_\theta}<\ln\left(\frac{K}{b(T)}\right)\}}\right]
}\\
&=&\PP\left({B_1}\geq \frac{1}{\sqrt\theta}\ln(\frac{K}{b(T)})\right)+\sigma\sqrt\theta\E\left({ B_1} 1_{\left\{{B_1}\geq\frac{\ln\frac{K}{b(T)}}{\sqrt\theta}\right\}}\right)+o(\sqrt\theta)
\\&\leq&\sqrt\theta \left(\frac{1}{\ln(\frac{K}{b(T)})}+\sigma\right)\frac{e^{-\frac{1}{2\theta}\ln^2(\frac{K}{b(T)})}}{\sqrt{2\pi}}+o(\sqrt\theta)
=o(\sqrt\theta)
.
\end{eqnarray*}
}\end{dema}
\begin{dema}\textbf{ of Lemma  \ref{min_deriv2_P}:\\}{
Let be $x\in \left(b(t),b(T)\right)$, then the variational inequality gives, for almost   $u\in \left(b(t),x\right)$,
\begin{eqnarray*}
\frac{u^2\sigma^2}{2}\frac{\partial^2P}{\partial x^2}(t,u)&\geq& rP(t,u)-(r-\delta)u\frac{\partial P}{\partial x}(t,u)\\
&&-\int \left (P(t,ue^y)-P(t,u)-u(e^y-1)\frac{\partial P}{\partial x}(t,u)\right)\nu(dy).\\
\end{eqnarray*}
Notice that $P(t,u)\geq K-u$, thus  
\begin{eqnarray}
\frac{u^2\sigma^2}{2}\frac{\partial^2P}{\partial x^2}(t,u)
&\geq&
r(K-u)+\left(r-\delta\right )u\nonumber\\
&&-\int \left (P(t,ue^y)-(K-u)+u(e^y-1)\right)\nu(dy)\nonumber\\
&&-u\left(\frac{\partial P}{\partial x}(t,u)+1\right)\left((r-\delta)-\int(e^y-1)\nu(dy)\right)\\
\end{eqnarray}
 And thanks to Lemma  \ref{deriv_prime},  we also have, for all $b(0)\leq u\leq x\leq b_e(t)\wedge b(T)$, $$\frac{\partial P}{\partial x}(t,u)+1=o(\sqrt{\theta}),$$ independently of $u$, therefore,
 \begin{eqnarray}
\frac{u^2\sigma^2}{2}\frac{\partial^2P}{\partial x^2}(t,u)
&\geq& rK-\delta u-\int \left (P(t,ue^y)-(K-ue^y)\right)\nu(dy)+o(\sqrt\theta). \label{a}
\end{eqnarray}
As the right hand side of  equality (\ref{a}) is non increasing in $u$, we obtain
\begin{eqnarray}
\inf_{b(t)\leq u \leq x}\frac{u^2\sigma^2}{2}\frac{\partial^2P}{\partial x^2}(t,u)
&\geq& rK-\delta x-\int \left (P(t,xe^y)-(K-xe^y)\right)\nu(dy)+o(\sqrt\theta).\nonumber\\\label{b}
\end{eqnarray}
Notice that \begin{eqnarray*}
{\int P(t,xe^y)\nu(dy)}&=&P_e(t,xe^y)+e(\theta,xe^y)\\
&=& \int \E(K-xe^ye^{X_\theta})^+\nu(dy)+o(\sqrt{\theta})\\
&=& \int \E\left (K-xe^y(1+\sigma B_\theta)\right )^+\nu(dy)+o(\sqrt{\theta})\\
&=& \int \E\left ((K-xe^y)-xe^y\sigma B_\theta\right )^+\nu(dy)+o(\sqrt{\theta}).\\
\end{eqnarray*}
We now consider  the integral $\int P(t,xe^y)\nu(dy)$
over  the sets $\{y<\ln(\frac{K}{b(T)})\}$,  $\{\ln(\frac{K}{b(T)})<y \}$ and $\{y=\ln(\frac{K}{b(T)})\}$.
Then, on the set $\{y<\ln(\frac{K}{b(T)}) \}$, we have  
\begin{eqnarray*}
\lefteqn{\int_{\{y<\ln(\frac{K}{b(T)}\}} P(t,xe^y)\nu(dy)
= \int_{\{y<\ln(\frac{K}{b(T)}\}} \hspace{-10mm}\E\left (K-xe^y)-xe^y\sigma B_\theta\right )^+\nu(dy)+o(\sqrt{\theta})}\\
&=& \int_{\{y<\ln(\frac{K}{b(T)}\}} (K-xe^y)\PP(xe^y\sigma B_\theta<(K-xe^y))\nu(dy)\\
&&-\int_{\{y<\ln(\frac{K}{b(T)}\}} xe^y\sigma \E\left (B_\theta1_{\{xe^y\sigma B_\theta<(K-xe^y)\}}\right )\nu(dy)+o(\sqrt{\theta})\\
&\leq& \int_{\{y<\ln(\frac{K}{b(T)}\}} \hspace{-10mm} (K-xe^y)
\nu(dy)-x\sigma\sqrt\theta\int_{\{y<\ln(\frac{K}{b(T)}\}}  \hspace{-10mm}e^y \E\left (B_1 1_{\{ B_1<\frac{1}{\sigma\sqrt{\theta}}(\frac{K}{x}e^{-y}-1)\}}\right )\nu(dy)+o(\sqrt{\theta}).\\
\end{eqnarray*}
For all $y<\ln(\frac{K}{b(T)})$, we have  $\frac{K}{x}e^{-y}-1>\frac{K}{b(T)}e^{-y}-1>0$, therefore
\begin{eqnarray*}
\lefteqn{0\leq -\E\left (B_1 1_{\{ B_1<\frac{1}{\sigma\sqrt{\theta}}(\frac{K}{x}e^{-y}-1)\}}\right )=\E\left (B_1 1_{\{ B_1\geq\frac{1}{\sigma\sqrt{\theta}}(\frac{K}{x}e^{-y}-1)\}}\right )}\\
&&\leq\E\left (B_1 1_{\{ B_1\geq\frac{1}{\sigma\sqrt{\theta}}(\frac{K}{b(T)}e^{-y}-1)\}}\right )
 \longrightarrow_{\theta \rightarrow0} 0.
\end{eqnarray*}
By the dominated convergence we obtain,
\begin{eqnarray}
{\int_{\{y<\ln(\frac{K}{b(T)}\}} P(t,xe^y)\nu(dy)}
&\leq& \int_{\{y<\ln(\frac{K}{b(T)}\}}(K-xe^y)\nu(dy)+o(\sqrt{\theta}).\label{set1}
\end{eqnarray}
On the set $\{y>\ln(\frac{K}{b(T)})\}$, we have $K<b(T)e^y$, therefore
\begin{eqnarray*}
\lefteqn{\int_{\{y>\ln(\frac{K}{b(T)}\}} \hspace{-5mm}P(t,xe^y)\nu(dy)= \int_{\{y>\ln(\frac{K}{b(T)}\}} \hspace{-5mm}\E\left ((K-xe^y)-xe^y\sigma B_\theta\right )^+\nu(dy)+o(\sqrt{\theta})}\\
&\leq&\int_{\{y>\ln(\frac{K}{b(T)}\}} \hspace{-3mm}\E\left[\left (b(T)e^y-xe^y-xe^y\sigma \sqrt\theta B_1 \right )1_{\{ xe^y\sigma\sqrt\theta B_1 <(K-xe^y) \}}\right]\nu(dy) +o(\sqrt{\theta})\\
&=& (b(T)-x)\int_{\{y>\ln(\frac{K}{b(T)}\}} e^y \PP \left( B_1<\frac{1}{\sigma \sqrt\theta}\left(\frac{K}{x}e^{-y}-1\right)\right)\nu(dy)\\
&&-\sqrt\theta x\int_{\{y>\ln(\frac{K}{b(T)}\}} e^y\sigma \E\left (B_1 1_{\{ B_1<\frac{1}{\sigma \sqrt\theta}\left(\frac{K}{x}e^{-y}-1\right)\}}\right )\nu(dy)+o(\sqrt{\theta})\\
\end{eqnarray*}
Notice that for all $y>\ln(\frac{K}{b(T)})$, we have $\frac{1}{\sigma \sqrt\theta}\left(\frac{K}{x}e^{-y}-1\right)\leq
\frac{1}{\sigma \sqrt\theta}(\frac{K}{b(t)}e^{-y}-1)\rightarrow-\infty$, thus
\begin{eqnarray*}
\PP(B_1<\frac{1}{\sigma \sqrt\theta}(\frac{K}{x}e^{-y}-1))&\leq&\PP(B_1<\frac{1}{\sigma \sqrt\theta}(\frac{K}{b(t)}e^{-y}-1))\\
&&\underset{\theta\rightarrow0}{\longrightarrow}0,
\end{eqnarray*}
and 
\begin{eqnarray*}
\E\left (|B_1| 1_{\{B_1<\frac{1}{\sigma \sqrt\theta}(\frac{K}{x}e^{-y}-1)\}}\right )
&\leq&\E\left (|B_1| 1_{\{B_1<\frac{1}{\sigma \sqrt\theta}(\frac{K}{b(t)}e^{-y}-1)\}}\right )\\
&&\underset{\theta\rightarrow0}{\longrightarrow}0.
\end{eqnarray*}
Therefore, by dominated convergence, we obtain
$$\int_{\{y>\ln(\frac{K}{b(T)}\}} e^y \PP(B_1<\frac{1}{\sigma \sqrt\theta}(\frac{K}{b(t)}e^{-y}-1))\nu(dy)\underset{\theta\rightarrow0}{\longrightarrow}0$$
and 
$$ -\sqrt\theta x\int_{\{y>\ln(\frac{K}{b(T)}\}} e^y\sigma \E\left (B_1 1_{\{xe^y\sigma \sqrt\theta B_1<(K-xe^y)\}}\right )\nu(dy)=o(\sqrt{\theta})$$
Consequently, if we denote by $\epsilon(\theta)=\int_{\{y>\ln(\frac{K}{b(T)}\}} e^y \PP(B_1<\frac{1}{\sigma \sqrt\theta}(\frac{K}{b(t)}e^{-y}-1))\nu(dy)$, we obtain
\begin{eqnarray}
\int_{\{y>\ln(\frac{K}{b(T)}\}} \hspace{-10mm}P(t,xe^y)\nu(dy)&\leq& (b(T)-x)\epsilon(\theta)
+o(\sqrt{\theta}),\label{set2}
\end{eqnarray}
with $\epsilon(\theta)\underset{\theta\rightarrow0}{\longrightarrow}0$.

Finally, on the set $\{y=\ln(\frac{K}{b(T)})\}$, 
we have \begin{eqnarray*}
\lefteqn{\int_{\{\ln(\frac{K}{b(T)})\}} P(t,xe^y)\nu(dy)= \int_{\{\ln(\frac{K}{b(T)})\}} \E\left ((K-xe^y)-xe^y\sigma B_\theta\right )^+\nu(dy)+o(\sqrt{\theta})}
\\&=& 
\int_{\{\ln(\frac{K}{b(T)})\}} \hspace{-8mm}(K-xe^y)\nu(dy)+\int_{\{\ln(\frac{K}{b(T)})\}} \hspace{-8mm}\E\left (xe^y\sigma B_\theta-(K-xe^y)\right )^+\nu(dy)+o(\sqrt{\theta})
\\&=&
   \int_{\{\ln(\frac{K}{b(T)})\}} \hspace{-8mm}(K-xe^y)\nu(dy)+ \int_{\{\ln(\frac{K}{b(T)})\}} \hspace{-8mm} xe^y\E\left (\sigma B_\theta-(\frac{K}{x}e^{-y}-1)\right )^+\nu(dy)
\\&=&
  \int_{\{\ln(\frac{K}{b(T)})\}} \hspace{-8mm}(K-xe^y)\nu(dy)+\frac{xK}{b(T)}\nu\left\{\ln \left(\frac{K}{b(T)}\right)\right\} \E\left (\sigma B_\theta-\left(\frac{b(T)}{x}-1\right)\right )^+
  \\&\leq&
  \int_{\{\ln(\frac{K}{b(T)})\}} \hspace{-8mm}(K-xe^y)\nu(dy)+K\nu\left\{\ln \left(\frac{K}{b(T)}\right)\right\} \E\left (\sigma B_\theta-\ln\left(\frac{b(T)}{x}\right)\right )^+
\end{eqnarray*}

We have thus proved that 
\begin{eqnarray*}
\int P(t,xe^y)\nu(dy)
&\leq&\int_{\{y\leq \ln(\frac{K}{b(T)}\}}\hspace{-12mm}(K-xe^y)\nu(dy)+K\nu\left\{\ln \left(\frac{K}{b(T)}\right)\right\} \E\left (\sigma B_\theta-\ln\left(\frac{b(T)}{x}\right)\right )^+
\\
&&+\left (b(T)-x\right )\epsilon (\theta)+o(\sqrt{\theta}),
\end{eqnarray*}
Coming back to  inequality (\ref{b}), we obtain
\begin{eqnarray*}
\lefteqn{
\int P(t,xe^y)-(K-xe^y)\nu(dy)
}\\
&\leq&-\int_{\{y> \ln(\frac{K}{b(T)}\}}\hspace{-12mm}(K-xe^y)\nu(dy)+K\nu\left\{\ln \left(\frac{K}{b(T)}\right)\right\} \E\left (\sigma B_\theta-\ln\left(\frac{b(T)}{x}\right)\right )^+
\\
&&+\left (b(T)-x\right )\epsilon (\theta)+o(\sqrt{\theta})
.
\end{eqnarray*}
Finally, since $rK=\delta b(T)+\int(b(t)e^y-K)^+\nu(y)$, we have
\begin{eqnarray*}
\lefteqn{\inf_{b(t)\leq u \leq x}\frac{u^2\sigma^2}{2}\frac{\partial^2P}{\partial x^2}(t,u)}\\
&\geq&rK-\delta x-\int \left (P(t,xe^y)-(K-xe^y)\right)\nu(dy)+o(\sqrt\theta)\\
&\geq& \hspace{-2mm}(b(T)-x)\hspace{-1mm}\left (\hspace{-1mm}\delta +\int_{\{y>\ln\frac{K}{b(T)}\}} \hspace{-14mm}e^y\nu(dy)+\epsilon (\theta)\hspace{-1mm}\right ) \hspace{-1mm}-\hspace{-1mm}K\nu\left\{\ln\frac{K}{b(T)}\right\}\E\left (\hspace{-1mm}\sigma B_\theta-\ln(\frac{b(T)}{x})\hspace{-1mm}\right )^+\hspace{-2mm}+o(\sqrt{\theta}).
\end{eqnarray*}
We note ${\alpha}=\frac{\nu\{\ln\left(\frac{K}{b(T)} \right)\}}{\bar\delta} \frac{K}{b(T)}$ and $\bar\delta=\delta +\int_{\{y>\ln(\frac{K}{b(T)}\}}  e^y\nu(dy),$
we then have  for all $u$ and all $x$ such that $b(t)\leq u\leq x<b(T)$ 
\begin{eqnarray*}
\lefteqn{\inf_{b(t)\leq u \leq x}\frac{u^2\sigma^2}{2}\frac{\partial^2P}{\partial x^2}(t,u)}\\
&\geq& b(T)\bar\delta \left (\frac{(b(T)-x)}{b(T)}-{\alpha} \E\left (\sigma B_\theta-\ln(\frac{b(T)}{x})\right )^+\right )-\left (b(T)-x \right )\epsilon (\theta)+o(\sqrt{\theta}).
\end{eqnarray*}
}
\end{dema}
\begin{remark}
The expression $\inf_{b(t)< u < x} \frac{u^2\sigma^2}{2}\frac{\partial^2P}{\partial x^2}(t,u)$ is justified thanks to the smoothness of $P$ in the continuation region which can be proved thanks to PDE arguments (see for instance \cite{smooth}). Nevertheless, we will only need this lower bound of the second derivative in the distribution sense ($\frac{\partial^2P}{\partial x^2}(t,du)$).  
\end{remark}

\section*{Appendix 2: A study of $v_{\lambda,\beta}$}
 \emph{\textbf{Lemma} \ref{v_alpha} : 
There exists $y_{\lambda,\beta}\in\left(0,( 1
                 +\lambda \beta(2+e^{\lambda})\right)$ such that  such that  
                 \[
 \forall y< -y_{\lambda,\beta},\quad
                   v_{\lambda,\beta}(y)=0.
 \] $$y_{\lambda,\beta}=-\inf\{x\in \R\;|\;v_{\lambda,\beta}(x)>0 \}. $$
}
\begin{dema}\textbf{ of Lemma \ref{v_alpha}:\\}
We have \[
v_{\lambda,\beta}(y)=  \sup_{\tau\in \T_{0,1}}\left(
    I_0(\tau)+I_1(\tau)\right),
\]
with 
\[
I_0(\tau)=\E\left( e^{\lambda\tau} 1_{\{\hat N_\tau=0\}}
\int_0^{\tau} f_{\lambda\beta}(
   y+ B_{s}
    )ds
    \right),
\]
and
\[
I_1(\tau)=\beta\E\left(e^{\lambda\tau} 1_{\{\hat N_\tau=1\}}
   \left((y+B_\tau)^+-(y+B_{\hat T_1})^+\right) \right).
\]
we will study $I_0(\tau)$ and  $I_1(\tau)$.
First of all, we note that the process $(M^0_t)_{t\geq 0}$
defined by $M^0_t=e^{\lambda t} 1_{\{\hat N_t=0\}}$ is a non negative martingale
with $M^0_0=1$. Under the probability
$\PP^0$ with  density  $M^0_t$ on $\F_t$, it is straightforward to check that 
$(B)_{t\geq 0}$ remains a  $\mathbb F$-Brownian motion. 
We have if $y\leq 0$,
\begin{eqnarray*}
I_0(\tau)&=&\E^0\left( 
           \int_0^{\tau} f_{\lambda\beta}(
   y+ B_{s}
    )ds
    \right)\\
    &=&\E^0 \left(
          y\tau+ \int_0^{\tau} B_s ds+{\lambda\beta}\int_0^{\tau}(
   y+ B_{s}
    )^+ds\right)
    \\
    &\leq&y \E^0\left( \tau\right)+
   (1+{\lambda\beta}) \E^0\left(
          \int_0^{\tau} B_{s}^+ds\right)\\
          &\leq&
       y \E^0\left(  \tau\right)+
   ({\lambda\beta}+1) \E^0\left( 
          \int_0^{\tau} \E^0\left(B_{\tau}^+\; |\; \F_s\right) ds\right).
\end{eqnarray*}
Notice that, for $\tau \in \T_{0,1}$,
\begin{eqnarray*}
 \E^0\left(
          \int_0^{\tau} \E^0\left(B_{\tau}^+\; |\; \F_s\right) ds\right)
          &=&
           \E^0\left( 
          \int_0^1 1_{\{\tau>s\}} \E^0\left(B_{\tau}^+\; |\; \F_s\right) ds\right)\\
          &=&
          \int_0^1\E^0\left( 1_{\{\tau>s\}}  
           \E^0\left(B_{\tau}^+\; |\; \F_s\right)\right)ds\\
                      &=&
           \E^0\left(\tau  
           B_{\tau}^+\right)\\
           &\leq &
           \E^0\left(\frac{\tau ^2+B_\tau^2}{2}\right)\leq \E^0(\tau),
\end{eqnarray*}
where, we used  $0\leq \tau\leq 1$, for the last inequality.
We then have
\begin{eqnarray*}
I_0(\tau)&\leq&(y+{\lambda\beta} +1) \E^0\left( \tau\right).
\end{eqnarray*}
For the study of  $I_1(\tau)$, let us introduce the martingale
$(M^1_t)_{0\leq t\leq 1}$ defined by
\begin{eqnarray*}
M^1_t&=& \E\left(e^{\lambda} 1_{\{\hat N_1=1\}} \;|\;\F_t\right)\\
       &=&
           \E\left(e^{\lambda} 1_{\{\hat N_1=1, \hat N_t=0\}} \;|\;\F_t\right)
            +\E\left(e^{\lambda} 1_{\{\hat N_1=1, \hat N_t=1\}} \;|\;\F_t\right)\\
            &=&
             1_{\{\hat N_t=0\}}e^\lambda \PP(\hat N_1-\hat N_t=1)
            + 1_{\{\hat N_t=1\}}e^\lambda \PP(\hat N_1-\hat N_t=0)\\
            &=&
      1_{\{\hat N_t=0\}}\lambda(1-t)e^{\lambda t}
            + 1_{\{\hat N_t=1\}}e^{\lambda t}.         
\end{eqnarray*}
Under the probability
$\PP^1$ with density  $M^1_t/\lambda$ on $\F_t$, 
 it is straightforward to check that 
$(B_t)_{0\leq t \leq 1}$ remains a  $\mathbb F$-Brownian motion. 
  We have for $y<0$,
\begin{eqnarray*}
I_1(\tau)&=&\lambda \beta\E^1\left( 
   (y+B_\tau)^+-(y+B_{\hat T_1\wedge \tau})^+\right) \\
       &\leq&
       \lambda \beta\E^1\left( 
   (y+B_\tau)^+\right)\\
   &\leq&
   \lambda \beta\E^1\left( 
   B_\tau 1_{\{B_\tau>-y\}}\right)\\
   &\leq&
   \lambda \beta\E^1\left( 
   B_\tau^2/|y|\right)=
        \lambda \beta\E^1\left( 
   \tau\right)/|y|.
\end{eqnarray*}
Using the two upper bound of   $I_0(\tau)$ and  $I_1(\tau)$,
we obtain 
\begin{eqnarray*}
v_{\lambda,\beta}(y)&\leq&
     \sup_{\tau\in \T_{0,1}}\left(
     (y+{\lambda\beta}+1)\E^0(\tau)
     +\frac{\lambda\beta}{|y|} \E^1(\tau)
               \right)\\
               &=&
               \sup_{\tau\in \T_{0,1}}\E\left(
     (y+{\lambda\beta}+1)\tau M^0(\tau)
     +\frac{\beta}{|y|} \tau M^1(\tau)
               \right)  \\
               &=&
                \sup_{\tau\in \T_{0,1}}\E\left(
     (y+{\lambda\beta}+1)\tau e^{\lambda \tau}  1_{\{\hat N_\tau=0\}}
     +\frac{\beta}{|y|} \tau 
     \left( 1_{\{\hat N_\tau=0\}}\lambda(1-\tau)e^{\lambda \tau}
            + 1_{\{\hat N_\tau=1\}}e^{\lambda \tau}\right)
               \right) 
               \\
               &\leq&
               \sup_{\tau\in \T_{0,1}}\E\left(f(\tau,\hat N_\tau)
               \right), 
\end{eqnarray*}
with
\[
f(t,x)= 1_{\{ x=0\}}te^{\lambda t}
    \left(y+1+\lambda\beta(1+\frac{1}{|y|})\right)
                          +
                           1_{\{ x=1\}}\beta te^{\lambda t}/|y|.
 \]
 Notice that
 \[
  \sup_{\tau\in \T_{0,1}}\E\left(f(\tau,\hat N_\tau)
               \right)=
                \sup_{\tau\in \T_{0,1}(\hat N)}\E\left(f(\tau,\hat N_\tau)
               \right),
 \]
where $ \T_{0,1}(\hat N)$ denotes the set of the stopping times of the natural completed filtration of the process  $(\hat N_t)_{t\geq 0}$,
 with values in $[0,1]$.
 
 Then, if $\tau\in \T_{0,1}(\hat N)$, there exists , thanks to Lemma~\ref{lem-Poisson}, $t_0\in [0,1]$,
such that  \[
  \tau\wedge \hat T_1=t_0\wedge \hat T_1.
  \]
  we then have
  \begin{eqnarray*}
\E\left(\tau e^{\lambda \tau}  1_{\{\hat N_\tau=0\}}\right)
      &=&\E\left(\tau e^{\lambda \tau}  1_{\{\hat T_1>\tau\}}\right)\\
      &=&t_0 e^{\lambda t_0}\PP(\hat T_1>\tau)\\
      &=&t_0,
\end{eqnarray*}
and
 \begin{eqnarray*}
\E\left(\tau e^{\lambda \tau}  1_{\{\hat N_\tau=1\}}\right)
      &\leq&\E\left(\tau e^{\lambda \tau}  1_{\{\hat T_1\leq \tau\}}\right)\\
      &=&\E\left(\tau e^{\lambda \tau}  1_{\{\hat T_1\leq t_0 \}}\right)\\
      &\leq&e^{\lambda}\PP(\hat T_1\leq t_0)\\
      &=&e^{\lambda}(1-e^{-\lambda t_0})\leq \lambda e^\lambda t_0.
\end{eqnarray*}
we deduce that
\begin{eqnarray*}
   \sup_{\tau\in \T_{0,1}(\hat N)}\E\left(f(\tau,\hat N_\tau)
               \right)&\leq&
               \sup_{0\leq t_0\leq 1}
               \left(t_0\left(y +1
                 +\lambda\frac{\beta(2+e^{\lambda})}{|y|}\right)\right).
\end{eqnarray*}
The right hand side of this equation will be equal to 0 if
\[
y+ 1  +\lambda\frac{\beta(2+e^{\lambda})}{|y|}\leq 0,
\]
and particularly, if $y\leq -\left(1
                 +\lambda \beta(2+e^{\lambda})\right)$, then$$-y_{\lambda,\beta}\geq -\left(1
                 +\lambda \beta(2+e^{\lambda})\right).$$
  
   To prove $-y_{\lambda,\beta}<0$, we consider $y=0$. Since for all stopping time $\tau$,
   \[
\E\left( e^{\lambda\tau} 1_{\{\hat N_\tau=0\}}
\int_0^{\tau} {\lambda\beta}(
   y+ B_{s}
    )^+ds
+\beta e^{\lambda\tau} 1_{\{\hat N_\tau=1\}}
   \left((y+B_\tau)^+-(y+B_{\hat T_1})^+\right) \right)\geq0.
\]
    we have $$\vv_{\lambda,\beta}(0)\geq \sup_{\tau \in \T_{0,1}}\E\int_0^\tau B_s ds=\vv_0(0),$$
and it is proved in \cite{villeneuve} or \cite{thesme}, Proposition 2.2.4. $\vv_0(0)>0$.
\end{dema}
\begin{lemma}\label{lem-Poisson}
   Let  $N=(N_t)_{t\geq 0}$ a homogenous  Poisson process
   with intensity  $\lambda$, and  $T_1$ its first jump time. if  $\tau $ is a stopping time of the natural completed filtration of $N$ such that
   $\tau\leq T_1$ a.s.,  then, $\tau=T_1$ a.s., or  
   there exists $t_0\geq 0$, such that 
   $\tau=t_0\wedge T_1$ a.s.
\end{lemma}
\begin{dema}
We denote by $\mathbb F=(\F_t)_{t\geq 0}$ the natural completed  filtration of $N$.
First of all, notice that for all $t\geq 0$ and  $A\in\F_t$,
\[
\PP(A\; | \; N_t=0)\in \{0,1\}.
\]
Indeed, the $A$ having this  property  form a sub $\sigma$-algebra of
$\F_t$ which contains the events of the form
$\{N_{s}=n\}$, with  
$0\leq s\leq t$ and $n\in\N$.
\\
Now, let  $\tau$ be a  $\mathbb F$-stopping time. We have for all
$t\geq 0$, $\PP(\tau>t\;|\;N_t=0)\in \{0,1\}$. We set
\[
I=\{t\in [0,+\infty[\; |\; \PP(\tau>t\;|\;N_t=0)=0\}.
\]
Notice that $t\in I$ if and only if $\PP(\tau >t, T_1>t)=0$, or 
\[
t\in I \Leftrightarrow \PP\left(\tau\wedge T_1\leq t\right)=1.
\]
If $\tau\leq T_1$ a.s. and if  $\PP(\tau<T_1)>0$, there exists $s>0$ (rational number)
such that $\PP(\tau\leq s, s<T_1)>0$, hence $\PP(\tau\leq s\;|\; N_s=0)>0$,
and  $\PP(\tau>s\;|\; N_s=0)=0$. We deduce that  $I$ is non empty and 
we can write
\[
I=[t_0,+\infty[, \mbox{ avec } t_0=\inf\{t\geq 0\;|\; \PP(\tau>t\;|\;N_t=0)=0\}.
\]
We then have $\tau\wedge T_1\leq t_0$ a.s., hence
$\tau\leq t_0\wedge T_1$. Moreover, for $s<t_0$,
we have  $\PP(\tau>s\;|\;N_s=0)=1$ and  $\PP(\tau \leq s\;|\;N_s=0)=0$,
hence $\PP(\tau\leq s, s<T_1)=0$.Therefore, $\PP(\tau<t_0\wedge T_1)=0$
and consequently $\tau= t_0\wedge T_1$ 
a.p.
\end{dema}

\emph{\textbf{Lemma} \ref{conjecture} :
For all  $x>y_{\lambda,\beta}$, we have  
$$C(x)>0.$$.}
  \begin{dema}
We have    $ v_{\lambda,\beta}(-y_{\lambda,\beta}) = 0$, considering the stopping time  $\tau=1$, we obtain
\[
\E\left[ e^{\lambda } 1_{\{\hat N_1=0\}}
\int_0^{1} f_{\lambda\beta}(
    B_{s}-y_{\lambda,\beta}
    )ds
     +{\beta} 
     e^{\lambda} 1_{\{\hat N_1=1\}}
      \left((B_1-y_{\lambda,\beta})^+)-(B_{\hat T_1}-y_{\lambda,\beta})^+
         \right) 
          \right]\leq 0.
    \]
    However, we have, using the independence  between $\hat N$ and  $B$, 
    \begin{eqnarray}
      \E\left[e^{\lambda }  1_{\{\hat N_1=0\}}
\int_0^{1} f_{\lambda\beta}(
    B_{s}-y_{\lambda,\beta}
    )ds
          \right]&=& e^{\lambda }\PP(\hat N_1=0)\left (-y_{\lambda,\beta} +{\lambda\beta}\E\int_0^1( B_{s}-y_{\lambda,\beta})^+ds\right )\nonumber\\
          &=&-y_{\lambda,\beta} +{\lambda\beta}\E\int_0^1( B_{s}-y_{\lambda,\beta})^+ds.\label{membr1ty}
    \end{eqnarray}
   On the other hand, we have
    \begin{eqnarray*}
    \lefteqn{\E\left[ {\beta} 
     e^{\lambda} 1_{\{\hat N_1=1\}}
      \left((B_1-y_{\lambda,\beta})^+)-(B_{\hat T_1}-y_{\lambda,\beta})^+
         \right) 
          \right]}
          \\&=& {\beta} e^{\lambda} \PP(\hat N_1=1) \left[   
     \E(B_1-y_{\lambda,\beta})^+-\E \left ((B_{\hat T_1}-y_{\lambda,\beta})^+    |    \hat N_1=1\right )
          \right]
                    \\&=& {\beta} {\lambda}  \left[   
     \E(B_1-y_{\lambda,\beta})^+-\E \left ((B_{\hat T_1}-y_{\lambda,\beta})^+    |    \hat T_1\leq 1\right )
          \right].
    \end{eqnarray*}
Noticing that  $\lambda \beta={\lambda\beta}$ and that conditionally to $\{ \hat T_1\leq 1\}$, $ \hat T_1$ is uniformly distributed on $[0,1]$, we obtain
 \begin{eqnarray}
    \lefteqn{\E\left[ {\beta} 
     e^{\lambda} 1_{\{\hat N_1=1\}}
      \left((B_1-y_{\lambda,\beta})^+)-(B_{\hat T_1}-y_{\lambda,\beta})^+
         \right) 
          \right]} \nonumber     
           \\&=& {\lambda\beta}  \left[   
     \E(B_1-y_{\lambda,\beta})^+-\E \left (\int_0^1 (B_{s}-y_{\lambda,\beta})^+ ds\right )
          \right].\label{membr2ty}
    \end{eqnarray}
    Combining  \eqref{membr1ty} and  \eqref{membr2ty}, we have
    \begin{equation*}
     -y_{\lambda,\beta}+{\lambda\beta}   \E(B_1-y_{\lambda,\beta})^+=-C(y_{\lambda,\beta})\leq 0
    \end{equation*}
    To conclude the proof, we  use the strict increasing of $C$, hence for all $x>y_{\lambda,\beta}$, we have $$C(x)>C(y_{\lambda,\beta})\geq0.$$
\end{dema}

\bibliographystyle
{amsplain}
\bibliography{bibthese}
\end{document}